\documentclass{article}
\usepackage[%
journal=JNCG,    
lang=british,   
]{ems-journal}

\usepackage{mathtools}
\usepackage[capitalize]{cleveref}

\usepackage{subcaption}
\usepackage{pgfplots}
\usetikzlibrary{patterns}

\usepackage{makecell} 

\pgfplotsset{compat=newest}

\theoremstyle{plain}  
\newtheorem{theorem}{Theorem}[section]

\newtheorem{lemma}[theorem]{Lemma}
\newtheorem{proposition}[theorem]{Proposition}
\newtheorem{corollary}[theorem]{Corollary}

\newtheorem{question}[theorem]{Question}

\theoremstyle{definition}  
\newtheorem{definition}[theorem]{Definition}
\newtheorem{example}[theorem]{Example}
\newtheorem{remark}[theorem]{Remark}

\theoremstyle{plain}
\newtheorem*{theoremA}{Theorem A}
\newtheorem*{theoremB}{Theorem B}

\newcommand{\WR}{\mathcal{WR}}

\DeclareMathOperator{\asdim}{asdim}
\DeclareMathOperator{\trasdim}{trasdim}
\DeclareMathOperator{\supp}{supp}
\DeclareMathOperator{\Ord}{Ord}
\DeclareMathOperator{\orb}{orb}

\numberwithin{equation}{section}

\begin{document}

\title{Asymptotic Property C for Certain Wreath-like Products of Groups}



\emsauthor{1}{
	\givenname{Han}
	\surname{Liu}
	\mrid{}
	\orcid{}}{H.~Liu}

\Emsaffil{1}{
        \department{School of Mathematical Sciences}
        \organisation{Fudan University}
        \rorid{013q1eq08}
        \address{220 Handan Road}
        \zip{200433}
        \city{Shanghai}
        \country{China}
        \affemail{h\_liu@fudan.edu.cn}
}

\classification[51F30,54F45]{20F69}


\begin{abstract}
In this paper, we present generalizations of some results on the asymptotic property C for wreath products. Specifically, we prove that certain wreath-like products admit asymptotic property C, thus providing some new examples for further study. 
\end{abstract}

\maketitle


\section{Introduction} 

In 1993, M. Gromov introduced the notion of asymptotic dimension for metric spaces (\textnormal{cf.} \cite{Gromov1991}) as a large-scale analog of the Lebesgue covering dimension. This concept later inspired the notion of finite decomposition complexity, introduced by E. Guentner, R. Tessera, and G. Yu in 2012, to study questions concerning the topological rigidity of manifolds (\textnormal{cf.} \cite{GTYu2012}). 
In 2000, A. Dranishnikov introduced asymptotic property C for metric spaces as a large-scale analog of topological property C (\textnormal{cf.} \cite{Dran2000}). 
It is well-known that every metric space with finite asymptotic dimension has asymptotic property C and finite decomposition complexity. 

In 2010, T. Radul generalized asymptotic dimension to transfinite asymptotic dimension, leading to an equivalent characterization of asymptotic property C (\textnormal{cf.} \cite{RT2010}). 

In 2023, J. Zhu and Y. Wu proved that $\mathbb{Z} \wr\mathbb{Z}$ has asymptotic property C by studying its transfinite asymptotic dimension, in particular proving that the transfinite asymptotic dimension of $\mathbb{Z} \wr \mathbb{Z}$ is at most $\omega + 1$ (\textnormal{cf.}~\cite{ZhuWu2023}). Consequently, $\mathbb{Z} \wr\mathbb{Z}$ became the first finitely generated group with asymptotic property C and infinite asymptotic dimension.

In this paper, we extend their results by providing more examples of groups with asymptotic property C, constructed via wreath-like products of groups. 
The content is organized as follows. In Section 2, we introduce basic definitions and properties related to asymptotic property C and wreath-like products. 
In Section 3, we propose a technical condition called the \emph{disjointness of the tiling through uniform translation} (DTUT; see Definition~\ref{def:(m,n)-DRUT}), then we present and prove the main theorem. In Section 4, we investigate groups satisfying the DTUT property. In Section 5, we discuss possible directions for future research and propose some open questions. 

\begin{theoremA}[Theorem \ref{thm: GwrlikeH has apc}]
Let $G$ and $H$ be countable discrete groups equipped with proper left-invariant metrics. Suppose that $G$ has finite asymptotic dimension and satisfies the $(m,0)$-DTUT property for some $m\in\mathbb{N}$. Suppose that $H$ is an FC group with finite asymptotic dimension.
Let $W$ be a wreath-like product of $G$ and $H$ corresponding to the action $H\curvearrowright I$ with finitely many orbits, as given by the following extension
\begin{align*}
1\longrightarrow \bigoplus_{i\in I}G\longrightarrow  W \stackrel{\varepsilon}\longrightarrow H\longrightarrow 1.
\end{align*}

Suppose $T$ is a symmetric transversal associated with $\bigoplus_{I}G\lhd W$ such that the transversal map $\eta: H\longrightarrow T$ is unital and bornologous. Equip $W$ with a proper left-invariant metric associated with $T$, then  $\trasdim(W) \le \omega {m} + \asdim(H)$, and hence $W$ has asymptotic property C. 
\end{theoremA}

The following theorem summarizes several concrete examples of groups that satisfy the DTUT property. Moreover, each of these groups is known to have finite asymptotic dimension (\textnormal{cf.}~\cite{Gromov1981,Osin2005,BD2008}).
\begin{theoremB}\label{main thm B}
The following groups have the DTUT property: 
\begin{enumerate}
    \item Every finitely generated torsion-free nilpotent group of class $m$ admits the $(m,0)$-DTUT property. More generally, every finitely generated group with polynomial growth has the $(m,0)$-DTUT property for some $m \in \mathbb{N}$ (\cref{thm:f.g. torsion-free group has DTUT} and Corollary~\ref{coro:gps with poly grow has DTUT}). 
    
    \item Every hyperbolic group has the $(1,0)$-DTUT property, and so does any group that is hyperbolic relative to a collection of subgroups, each of which quasi-isometrically embeds into a finite product of metric trees (Proposition~\ref{prop:F_2 has (1,0)-DTUT} and Corollary~\ref{coro:hyper gps has (1,0)-DTUT}).
\end{enumerate}
\end{theoremB}

\section{Preliminaries}
\subsection{Asymptotic Property C}
Let $(X, d)$ be a metric space. For $U,V\subseteq X$, define 
\[
\text{diam}~ U\coloneqq \sup\{d(x,y) \mid  x,y\in U\}, \quad
d(U,V)\coloneqq \inf\{d(x,y) \mid  x\in U,y\in V\}.
\]
Let $R>0$ and $\mathcal{U}$ be a family of subsets of $X$, $\mathcal{U}$ is said to be \emph{$R$-bounded} if
\[
\text{diam}~\mathcal{U}\coloneqq\text{sup}\,\{\text{diam}~ U \mid  U\in \mathcal{U}\}\le R.
\]
The family $\mathcal{U}$ is said to be \emph{uniformly bounded} if there exists $R>0$ such that $\mathcal{U}$ is $R$-bounded. For $r>0$, $\mathcal{U}$ is said to be\emph{ $r$-disjoint} if $d(U,V)\ge r$ for all $U,V\in \mathcal{U}$ with $U\neq V$.

Throughout this paper, we use the following notations:
\begin{center}
\begin{tabular}{@{} >{\raggedright\arraybackslash}p{0.25\textwidth} p{0.70\textwidth} @{}}
\toprule
\textbf{Notation} & \textbf{Meaning} \\
\midrule
$\bigcup\mathcal{U}$ &  The union of all sets in the family $\mathcal{U}$, i.e., $\bigcup_{U \in \mathcal{U}} U$ \\
$\mathcal{U}_{1}\cup\mathcal{U}_{2}$ & The union of two families, i.e., $\{ U \mid U \in \mathcal{U}_1 \text{ or } U \in \mathcal{U}_2 \}$ \\
$\bigcup_{i\in\mathbb{N}^+}\mathcal{U}_{i}$ & The union of countably many families, \\
 &i.e., $\{ U \mid U \in \mathcal{U}_i \text{ for some } i \in \mathbb{N}^+ \}$ \\
$\mathcal{U}\setminus V$ & The family of set differences, i.e., $\{ U \setminus V \mid U \in \mathcal{U} \}$\\
$\mathcal{U} \cap V$ & The family of intersections, i.e., $\{ U \cap V \mid U \in \mathcal{U} \}$ \\
\bottomrule
\end{tabular}
\end{center}

Moreover, given a group $\Gamma$ acting on $X$, we write $a{.}x$ to denote the action of $a \in \Gamma$ on $x \in X$, and $a \cdot b$ to denote the group multiplication in $\Gamma$ (the symbol ``$\cdot$'' may sometimes be omitted). When $\mathcal{U}$ is a family of subsets of $\Gamma$, we naturally write $a \cdot \mathcal{U}$ to denote $\{a \cdot U \mid U \in \mathcal{U}\}$.

\begin{definition}[finite asymptotic dimension, \textnormal{cf.} \cite{Gromov1991}]
A metric space $X$ is said to have \emph{finite asymptotic dimension} if there exists $n\in\mathbb{N}$, such that for every $r>0$, there exists a sequence of uniformly bounded families $\{\mathcal{U}_{i}\}_{i=0}^{n}$ of subsets of $X$ such that the family
$\bigcup_{i=0}^{n}\mathcal{U}_{i}$ covers $X$ and $\mathcal{U}_{i}$ is $r$-disjoint for all $0\le i\le n$. In this case, we say that $\asdim(X)\le n$. Moreover, we say that $\asdim(X) = n$ if $n$ is the smallest such integer for which $X$ admits the above decomposition.
\end{definition}

\begin{definition}[asymptotic property C, \textnormal{cf.} \cite{Dran2000}]
A metric space $X$ is said to have \emph{asymptotic property C} if for every infinite sequence $R_0<R_1<\cdots$ of positive numbers, there exist $n\in\mathbb{N}$ and uniformly bounded families $\mathcal{U}_0,\cdots,\mathcal{U}_n$ of subsets of $X$ such that the family $\bigcup_{i=0}^n\mathcal{U}_i$ covers $X$ and $\mathcal{U}_i$ is $R_i$-disjoint for all $0\le i\le n$.
\end{definition}

\begin{definition}[bornologous map, \textnormal{cf.} \cite{Willett2009}]
Let $f : X\longrightarrow Y$ be a map of metric spaces, we call $f$ \emph{bornologous} if for all $R > 0$ there exists $P(R) > 0$ such that if $d_X(x_1, x_2)<R$, then $d_Y(f(x_1),f(x_2)) < P(R)$. 
\end{definition}

Without loss of generality, we always assume that $P(R)$ is monotonically increasing. 

\begin{definition}[coarse equivalence, \textnormal{cf.} \cite{NY}]
Let $f : X\longrightarrow Y$ be a map between metric spaces, we call $f$ a \emph{coarse equivalence} if 
\begin{enumerate}
    \item there exist non-decreasing functions $\rho_+,\rho_-:[0,\infty]\longrightarrow [0,\infty]$   such that 
     \[
     \lim_{t\to\infty}\rho_-(t)=\infty,
     \]
     and the inequality
     \[
     \rho_-( d_X(x, y))\le d_Y(f(x),f(y))\le  \rho_+( d_X(x, y)) 
     \]
     holds for all $x,y\in X$, and 
    \item the image $f(X)$ is a net in $Y$ .
\end{enumerate}
Two metric spaces $X$ and $Y$ are called \emph{coarsely equivalent} if there exists a coarse equivalence $f:X\longrightarrow Y$. 
A map which satisfies (1), but not necessarily (2), is called a \emph{coarse embedding} of $X$ into $Y$.
\end{definition}

Note that a coarse embedding of $X$ into $Y$ is the same as a coarse equivalence of $X$ with a subspace of $Y$. All the large-scale properties mentioned in this paper (including finite asymptotic dimension and asymptotic property C) are coarse invariants (\textnormal{cf.} \cite{NY}).

\subsection{Transfinite Asymptotic Dimension}

In 2010, T. Radul generalized asymptotic dimension of a metric space $X$ to transfinite asymptotic dimension which is denoted by $\trasdim(X)$ (\textnormal{cf.} \cite{RT2010}).

Let $\text{Fin}~\mathbb{N}$ denote the collection of all finite, nonempty subsets of $\mathbb{N}$ and let $ M \subset \text{Fin}~\mathbb{N}$. For $\sigma\in \{\varnothing\}\bigcup \text{Fin}~\mathbb{N}$, let
\[
M^{\sigma} \coloneqq \{\tau\in \text{Fin}~\mathbb{N}  \mid  \tau \cup \sigma \in M \text{ and } \tau \cap \sigma = \varnothing\}.
\]

Let $M^a$ denote an abbreviation of $M^{\{a\}}$ for $a \in \mathbb{N}$. Define ordinal number $\Ord M$ inductively as follows:
\begin{align*}
\Ord M = 0 &\Leftrightarrow M = \varnothing,\\
\Ord M \le \alpha &\Leftrightarrow \forall~ a\in \mathbb{N}, ~\Ord M^a < \alpha,\\
\Ord M = \alpha &\Leftrightarrow \Ord M \le \alpha \text{ and } \Ord M < \alpha \text{ is not true},\\
\Ord M = \infty &\Leftrightarrow \Ord M \le\alpha \text{ is not true for every ordinal number } \alpha.
\end{align*}

Given a metric space $(X, d)$, define the following collection:
\[
\begin{split}
A(X, d) \coloneqq \{\sigma \in \text{Fin}~\mathbb{N}  \mid &\text{ there are no uniformly bounded families } \mathcal{U}_i  \text{ for } i \in \sigma\\
& \text{ such that each } \mathcal{U}_i \text{ is } i\text{-disjoint and }\bigcup_{i\in\sigma}\mathcal{U}_i \text{~covers~} X\}.
\end{split}
\]

\begin{definition}[transfinite asymptotic dimension, \textnormal{cf.} \cite{RT2010}]
The \emph{transfinite asymptotic dimension} of $X$ is defined as $\trasdim (X)\coloneqq \Ord \,A(X, d)$.
\end{definition}

\begin{remark}
Note that transfinite asymptotic dimension is a generalization of finite asymptotic dimension.
That is, $\trasdim(X)\le n$ if and only if $\asdim(X)\le n$ for each $n\in\mathbb{N}$.
\end{remark}

\begin{lemma}[\textnormal{cf.} \cite{RT2010}]
\label{lemma:equi between apc and trandim}
Let $X$ be a metric space, $X$ has asymptotic property $C$ if and only if $\trasdim(X)\le\alpha$ for some countable ordinal number $\alpha$.
\end{lemma}

Let $\omega$ be the smallest infinite ordinal. 
\begin{lemma}[\textnormal{cf.} \cite{WuZhu2018}]
\label{lemma: decomp for trasdim}
Let $X$ be a metric space with $\asdim(X) = \infty$ and $k\in \mathbb{N}$, the following are equivalent:
\begin{itemize}
\item $\trasdim(X) \le \omega+k$;
\item For every $n\in\mathbb{N}$, there exists $m(n)\in\mathbb{N},$ such that for every~ $d>0$, there are uniformly bounded families $\mathcal{U}_{-k},\mathcal{U}_{-k+1},\cdots,\mathcal{U}_{m(n)}$ such that $\mathcal{U}_i$ is $n$-disjoint for $-k\le i\le 0$, $\mathcal{U}_j$ is $d$-disjoint for $1\le j\le m(n)$ and
$\bigcup_{i=-k}^{m(n)}\mathcal{U}_i$ covers $X$. Moreover, $m(n)\rightarrow \infty$ as $n\rightarrow \infty$.
\end{itemize}
\end{lemma}

\subsection{Saturated Union}
\begin{definition}[\textnormal{cf.} \cite{BD2011}]
Let $\mathcal{U}$ and $\mathcal{V}$ be families of subsets of $X$. The \emph{$r$-saturated union} of $\mathcal{V}$ with $\mathcal{U}$ is defined as 
\[
\mathcal{V}\cup_r\mathcal{U}\coloneqq\{N_r(V;\mathcal{U}) \mid V\in\mathcal{V}\}\cup\{U\in\mathcal{U} \mid d(U,\mathcal{V})>r\},
\]
where $N_r(V;\mathcal{U}) \coloneqq V \cup\,\bigcup_{d(U,V)\le r} U$ and $d(U,\mathcal{V})>r$ means that $d(U,V)>r$ for all $V\in\mathcal{V}$.
\end{definition}

\begin{lemma}[\textnormal{cf.} \cite{BD2011}]
\label{lemma: r-saturated}
Let $\mathcal{U}$ be an $r$-disjoint and $R$-bounded family of subsets of $X$ with $R\ge r$. Let $\mathcal{V}$ be a $5R$-disjoint, $D$-bounded family of subsets of $X$. Then the family $\mathcal{V}\cup_r\mathcal{U}$ is $r$-disjoint, $(D+2R+2r)$-bounded and $\mathcal{V}\cup_r\mathcal{U}$
covers $\bigcup(\mathcal{V}\cup\mathcal{U})$.
\end{lemma}

\subsection{Groups as metric spaces}

Recall that a metric on a set is said to be \emph{proper} if closed balls are compact. A metric on a group is called \emph{left-invariant} if the action of the group on itself by left multiplication is an isometry. Given a group $G$ with a generating set $S$, let $1_{G}\in G$ be the unit element of $G$, for any $ g\in G$, let $l_G(g)$ to be the length of the shortest word representing $g$ in elements of $S\cup S^{-1}$. We say that $l_G$ is a \emph{word length function} for $G$. 
The left-invariant word metric $d_{G}$ on $G$ is induced by word length function, i.e.,
for every $ g, h\in G$,
\[
d_{G}(g, h)\coloneqq l_G(g^{-1}h).
\]

Finitely generated groups carry a unique (up to coarse equivalence) left-invariant proper metric. Considering a countable (discrete) group $G$ that is not finitely generated, A. Dranishnikov and J. Smith showed that, up to coarse equivalence, $G$ admits a unique proper left-invariant metric (\textnormal{cf.} \cite{DS2006}). This metric is constructed by taking a (countably) infinite symmetric generating set 
$S$ and defining a weighted word metric, where the weight function associated with the elements of $S$ is proper, i.e., for any $N\in \mathbb{N}$ the set of $s\in S$ with weight at most $N$ is finite. 

Thus, when we focus on proper left-invariant metrics, the coarse geometric properties of any countable group become intrinsic properties of the group itself. In this paper, we assume that all groups are equipped with proper left-invariant metrics unless otherwise stated.

Recall some permanence properties of groups with finite asymptotic dimension.
\begin{theorem}[\textnormal{cf.} \cite{BD2001,BD2006}]
The collection of countable groups having finite asymptotic dimension is closed under operations such as taking subgroups, direct products, extensions, free amalgamated products and HNN-extensions.	
\end{theorem}

Unlike groups with finite asymptotic dimension, permanence properties for asymptotic property C are more elusive.
\begin{theorem}[\textnormal{cf.} \cite{BN2018}]\label{thm: closed properties for APC}
The collection of countable groups having asymptotic property C is closed under operations such as taking subgroups, direct products and free products. Moreover, let $G$,$H$ and $K$ be groups with proper left-invariant metrics, if $1\to K\to G\to H\to 1$ is an exact sequence, $H$ has asymptotic property C, and $K$ has finite asymptotic dimension, then $G$ has asymptotic property C.
\end{theorem}

\subsection{Wreath-like Products of Groups}

Consider a group extension of $N$ by $H$ of the form $1\longrightarrow N \stackrel{\mu}\longrightarrow  W \stackrel{\varepsilon}\longrightarrow H\longrightarrow 1$. To simplify the notation without causing any misunderstanding, we will always omit $\mu$ and treat $N$ directly as a normal subgroup of $W$. 

In 2023, I. Chifan, A. Ioana, D. Osin and B. Sun introduced a new class of groups called wreath-like products, which are closely related to the classical wreath products and arise naturally in the context of group theoretic Dehn filling.
\begin{definition}[\textnormal{cf.} \cite{ClOS2023}]\label{def: wlp}
Let $G$ and $H$ be groups, and let $I$ be a set equipped with a (left) action of $H$ on $I$.
We say that a group $W$ is a \emph{wreath-like product} of groups $G$ and $H$ corresponding to the action $H\curvearrowright I$ if $W$ is an extension of the form
\begin{equation}\label{ext}
1\longrightarrow \bigoplus_{i\in I}G_i \longrightarrow  W \stackrel{\varepsilon}\longrightarrow H\longrightarrow 1, 
\end{equation}
where $G_i\cong G$ and the action of $W$ on $G^{(I)}=\bigoplus_{i\in I}G_i$ by conjugation satisfies the rule 
\begin{equation*}
wG_iw^{-1} = G_{\varepsilon(w){.}i}
\end{equation*}
for all $i\in I$. The map $\varepsilon\colon W\to H$ is called the \emph{canonical homomorphism associated with the wreath-like structure of $W$}. 

If the action $H\curvearrowright I$ is regular (i.e., free and transitive), we say that $W$ is a \emph{regular wreath-like product} of $G$ and $H$. The set of all wreath-like products of $G$ and $H$ corresponding to an action $H\curvearrowright I$ (respectively, all regular wreath-like products) is denoted by $\WR(G, H\curvearrowright I)$ (respectively, $\WR(G,H)$).
\end{definition}

The notion of a wreath-like product generalizes the ordinary (restricted) wreath product of groups. Indeed, given groups $G$ and $H$, and an action $ H \curvearrowright I $, the restricted wreath product is denoted by $ G\,\text{wr}_I\,H $. When $I = H$ and $H \curvearrowright H$ is the natural group action, we write it simply as $G \wr H$. Obviously, $G\,\text{wr}_I\, H\in \WR(G, H \curvearrowright I)$. Conversely, it is not difficult to show that $W\cong G\,\text{wr}_I\, H$ whenever the extension (\ref{ext}) splits.

Now we introduce the group structure of a group extension. 
Given a subgroup $K$ of a group $W$, a right (respectively, left) \textit{transversal} is a set containing exactly one element from each right (respectively, left) coset of $K$. In this case, the ``sets'' (cosets) are mutually disjoint, i.e., the cosets form a partition of the group $W$. In the following, we assume that every transversal refers to a right transversal. 

Given a group extension $1\longrightarrow N \longrightarrow  W \stackrel{\varepsilon}\longrightarrow H\longrightarrow 1$, we can express every element of $W$, uniquely up to the choice of a transversal, as a pair in the Cartesian product $N \times H$. 
Given a transversal $T$ associated with $N \lhd W$, and for each $w \in W$, there exist unique elements $n \in N$ and $t \in T$ such that $w = nt$. Since $W / N \cong H$, this determines a bijection $\eta: H \to T$, which we call a \emph{transversal map}. 
If moreover $\eta$ maps the unit element of $H$ to the unit element of $W$, we say that $\eta$ is \emph{unital}. 
Recall that a subset $T$ of a group $\Gamma$ is called \emph{symmetric} if $1_\Gamma \in T$ and $T$ is closed under taking inverses, i.e., $a \in T$ implies $a^{-1} \in T$.
It is worth noting that one can always choose a transversal $T$ that is symmetric, and the corresponding transversal map $\eta$ can be taken to be unital.

Thus, $w = n\eta(\eta^{-1}(t))$, and hence $w$ can be uniquely represented as a pair $(n, \eta^{-1}(t))\in N\times H$, up to the choice of the transversal $T$. It remains to describe the group multiplication in $W$ in terms of the Cartesian product $N \times H$ under this identification.

Let $w_1 = (n_1, h_1)$ and $w_2 = (n_2, h_2)$ be two elements in $W$, then we have
\begin{align*}
w_1w_2=&n_1\eta(h_1)\cdot n_2\eta(h_2)=n_1(\eta(h_1)n_2(\eta(h_1))^{-1})\eta(h_1)\eta(h_2)\\
=&n_1\big(\eta(h_1)n_2(\eta(h_1))^{-1}\big)\eta(h_1)\eta(h_2)\big(\eta(h_1h_2)\big)^{-1}\eta(h_1h_2).
\end{align*}
Note that $\eta(h_1)n_2(\eta(h_1))^{-1}\in N$ and $\eta(h_1)\eta(h_2)(\eta(h_1h_2))^{-1}\in N$, thus 
\begin{align*}
w_1w_2=\bigg(n_1\big(\eta(h_1)n_2(\eta(h_1))^{-1}\big)\eta(h_1)\eta(h_2)(\eta(h_1h_2))^{-1},~h_1h_2\bigg).
\end{align*}

For $W\in \WR(G, H\curvearrowright I)$ with a symmetric transversal $T\subset W$ and a unital transversal map $\eta: H\longrightarrow T$, note that for all $h\in H$, $\varepsilon\eta(h)=h$, and $\eta(h^{-1})=(\eta(h))^{-1}$.
Let $w_1=((a_i)_{i\in I},h_1)$ and $w_2=((b_i)_{i\in I},h_2)$ be elements of $W$, then the multiplication in $W$ is given by  
\begin{align*}
w_1w_2=\bigg((a_i)_{i\in I}(b_{h_1{.}i})_{i\in I}\eta(h_1)\eta(h_2)(\eta(h_1h_2))^{-1}, ~h_1h_2\bigg),
\end{align*}
and for all $w=((a_i)_{i\in I},h)\in W$, the inverse is given by
\begin{align*}
w^{-1}=\big( [(a_{h^{-1}{.} i})^{-1}]_{i\in I}, ~h^{-1}\big).
\end{align*}

Now we consider the generating set of $W$. 
Let $1_{G}\in G$ and $1_{H}\in H$ denote the unit elements of $G$ and $H$, respectively. Let $S_G$ and $S_H$ be generating sets of $G$ and $H$, respectively. Let $e \in \bigoplus_{i \in I} G_i$ denote the element that takes the value $1_G$ at every $i \in I$, i.e., $e$ is the unit element of $\bigoplus_{i\in I}G_i$. For every $i\in I$ and $b\in G$, let $\delta^{b}_{i}\in \bigoplus_{i\in I}G_i$ denote the $\delta$-function, defined by 
\[
\delta^{b}_{i}(i)\coloneqq b
\text{ and }\delta^{b}_{i}(j) \coloneqq 1_{G} \text{ for all }j\neq i.
\]

For each $w=((a_i)_{i\in I},h)\in W$, suppose $\supp((a_i)_{i\in I})=\{i_1,i_2,\cdots,i_n\}\subset I$, then 
\begin{align*}
w=&((a_i)_{i\in I},h)=((a_i)_{i\in I},1_H)\cdot (e,h)\\
=&(\delta^{a_{i_1}}_{i_1},1_H)(\delta^{a_{i_2}}_{i_2},1_H)\cdots (\delta^{a_{i_n}}_{i_n},1_H)\cdot (e,h)
\end{align*}

Denote the $H$-orbits in $I$ by $\{I_\alpha\}_{\alpha\in \orb(I)}$, where $\orb(I)$ is the index set. Throughout this paper, we always fix an element $i_{\alpha,0}\in I_\alpha$ for each $\alpha\in \orb(I)$. Then for any $i\in I$, there exist $\alpha\in \orb(I)$ and $p\in H$ such that $p{.} i_{\alpha,0}=i$. Moreover for any $g=s_1s_2\cdots s_n\in G$, we have 
\begin{align*}
(\delta^{g}_{i}, 1_H)=(e,p)(\delta^{s_1}_{i_{\alpha,0}},1_H)(\delta^{s_2}_{i_{\alpha,0}},1_H)\cdots(\delta^{s_n}_{i_{\alpha,0}},1_H)(e,p^{-1}).
\end{align*}

Let $\widetilde{S}=\big\{(\delta^{s_g}_{i_{\alpha,0}},1_H),(e,h) \mid  s_g\in S_G, h\in H, \alpha\in \orb(I)\big\}$, then $\widetilde{S}$ clearly generates $W$. 
Assign weight $l_H(h)$ to each generator of the form $(e,h)$ for all $h \in H$, and assign weight $l_G(s_g)$ to each generator of the form $(\delta^{s_g}_{i_{\alpha,0}},1_H)$ for all $s_g\in S_G$ and $\alpha\in \orb(I)$. Then the word metric $d_W$ associated with this weighted generating set defines a proper left-invariant metric on $W$, which we refer to as a metric associated with the transversal $T$. 

\begin{remark}\label{remark: generating set of wr product}
When both $G$ and $H$ are finitely generated, and $W$ is the classical restricted wreath product with finite $\orb(I)$, then $W$ can be finitely generated. A finite generating set can be chosen as
\[
{S}=\big\{(\delta^{s_g}_{i_{\alpha,0}},1_H),(e,s_h) \mid  s_g\in S_G, s_h\in S_H, \alpha\in \orb(I)\big\}. 
\]
In this case, for any $w=\big((a_i)_{i\in I},h\big)\in W$, we have $l_S(w)=l_{W, \widetilde{S}}(w)$.

Furthermore, for each $j\in \supp(a_i)_{i\in I}$, there exists a unique $\alpha\in \orb(I)$ such that $j\in I_\alpha$. Denote
$\widetilde{h_j^w}=\min\{h_j\in H \mid  h_j{.}i_{\alpha,0}=j \}$. Let $L_H(w)$ denote the length of the shortest path in the Cayley graph of $H$ starting from $1_H$, passing through all elements $\widetilde{h_j^w}$ with $j\in \supp(a_i)_{i\in I}$, and finally ending at $h$. Then we have
\begin{align*}
l_W(w)=\sum_{j\in \supp(a_i)_{i\in I}}l_G(a_j)+L_H(w). 
\end{align*}

Note that when $W$ is the classical restricted wreath product, for any $h_1, h_2 \in H$, we have $(e,h_1)\cdot (e,h_2) = (e,h_1h_2)$. However, this equality generally fails for wreath-like products that are not classical. Instead, the product yields an additional element $\eta(h_1)\eta(h_2)\eta(h_1h_2)^{-1}$ in $\bigoplus_{i\in I}G$, which may shorten the minimal word length of elements in the wreath-like product.   
\end{remark}

\begin{remark}
Note that for all $w_1=((a_i)_{i\in I},h), w_2=((b_i)_{i\in I},p)\in W$, we have
\begin{align*}
d_W(w_1,w_2)=&l_W\bigg( \big([(a_{h^{-1}{.} i})^{-1}]_{i\in I}, h^{-1}\big)\cdot \big((b_i)_{i\in I},p\big) \bigg)\\
=&l_W\bigg( \big[(a_{h^{-1}{.} i})^{-1}\cdot b_{h^{-1}{.} i}\big]_{i\in I}\cdot \eta(h^{-1})\eta(p)\eta(h^{-1}p)^{-1},~ h^{-1}p \bigg)
\end{align*}
\end{remark}

\begin{remark}\label{remark:bornologous}
The transversal map $\eta:(H,d_H)\longrightarrow (W,d_W)$ is bornologous if and only if for all $R>0$, there exist $P(R)>0$ and finitely many indices $\alpha_1,\alpha_2,\cdots,\alpha_{P(R)}\in \orb(I)$ such that for all $x,y\in H$, whenever $d_H(x,y)<R$, we have 
\begin{align*}
    \supp(\eta(x^{-1})\eta(y)\eta(x^{-1}y)^{-1})\subset \bigcup_{1\le j\le P(R)}B_H(1_H,P(R)){.} i_{\alpha_j,0}.
\end{align*} 
Moreover, for all $i\in \supp(\eta(x^{-1})\eta(y)\eta(x^{-1}y)^{-1})$, it holds that
\begin{align*}
l_{G}\bigg(\Big(\eta(x^{-1})\eta(y)\eta(x^{-1}y)^{-1}\Big)_i\bigg)\le P(R).
\end{align*}

The function $P:\mathbb{R}\to\mathbb{R}$ is called \emph{a control function associated with $\eta$}.

In addition, since for all $h_1,h_2\in H$, we have $d_W\big( (e,h_1),\,(e,h_2) \big)\ge d_H(h_1,h_2)$, it follows that $\eta$ is bornologous if and only if $\eta$ is a coarse embedding. 
\end{remark}

According to Remark \ref{remark: generating set of wr product}, for any classical wreath product $W = G\,\text{wr}_I\, H$, there always exists an isometric transversal map. In this paper, we weaken this requirement in the context of wreath-like products by requiring the transversal map to be a coarse embedding.

\section{Technical Tools and the Main Result}

Note that finite decomposition complexity (\textnormal{cf.} \cite{GTYu2012}) and asymptotic property C are both weaker than finite asymptotic dimension; however, the relationship between these two properties remains unknown. Currently, no example is known that distinguishes asymptotic property C from finite decomposition complexity.

Let $G$ and $H$ be groups with finite decomposition complexity, and let $W\in \WR(G, H\curvearrowright I)$, where $I$ is a countable set. It is known that the class of groups with finite decomposition complexity is closed under operations such as taking subgroups, extensions, free amalgamated products, HNN extensions, and direct unions (\textnormal{cf.} \cite{GTYu2013}). Therefore, $W$ inherits finite decomposition complexity. This naturally raises the question: under what conditions on $G$ and $H$ does $W$ have asymptotic property C? 

On the other hand, G. Bell, D. G\l odkowski, and A. Nag\'{o}rko raised the question (\textnormal{cf.} \cite{BGN2019}) of whether the class of countable groups with left-invariant proper metrics and asymptotic property C is closed under taking extensions and wreath products. 

In the following, we aim to provide new examples of groups with asymptotic property C that are constructed via wreath-like products, and to further explore the above questions concerning permanence properties of asymptotic property C.

We begin by extending \cref{lemma: decomp for trasdim} to the following more general case:
\begin{lemma}\label{lemma:more general decomp for trasdim}
Let $(X,d)$ be a metric space and $m, n\in \mathbb{N}$, the following are equivalent:

\begin{enumerate}
\item $\trasdim(X) \le \omega {m}+n$;
\item For every sequence $k_0, k_1, \cdots, k_m \in \mathbb{N}^+$ with $k_0 < k_1 < \cdots < k_m$, there exist integers $E_i = E_i(k_0, k_1, \cdots, k_{i-1}) \in \mathbb{N}$ for each $1 \le i \le m$, and uniformly bounded families
\[
\mathcal{U}_{0}, \mathcal{U}_{1}, \cdots, \mathcal{U}_{n},  \quad \mathcal{U}_{i,0}, \mathcal{U}_{i,1}, \cdots, \mathcal{U}_{i,E_i}  \quad (1 \le i \le m)
\]
satisfying the following: 
  \begin{itemize}
      \item Each $\mathcal{U}_j$ is $k_0$-disjoint for all $0\le j \le n$.
      \item Each $\mathcal{U}_{i,p_i}$ is $k_i$-disjoint for all $1 \le i \le m$ and $0 \le p_i \le E_i$.
      \item The family $\big(\bigcup_{j=0}^{n}\mathcal{U}_j\big) \cup \big(\bigcup_{i=1}^m\bigcup_{p_i=0}^{E_i}\mathcal{U}_{i,p_i}\big)$ covers $X$. 
  \end{itemize}
\end{enumerate} 
\end{lemma}

\begin{proof}
According to the definition of $\trasdim (X)$, we have 
\begin{align*}
\trasdim(X)\le \omega {m}+n 
\iff& \forall~ a^0_0<a^0_1<\cdots<a^0_n\in\mathbb{N}^+, \exists~ M_1=M_1(a_n^0)\in\mathbb{N}^+, s.t. \\
& \Ord A(X,d)^{\big\{ a^0_{0},a^0_{1},\cdots,a^0_{n} \big\}}\le \omega {(m-1)} + M_1.
\end{align*}
Then, by induction, we further obtain 
\begin{align*}
&\trasdim(X)\le \omega {m}+n \\
\iff& \forall~ a^0_0<a^0_1<\cdots<a^0_n\in\mathbb{N}^+, \exists~ M_1=M_1(a_n^0)\in\mathbb{N}^+, s.t. \\
& \quad\forall~ a^1_0<a^1_1<\cdots<a^1_{M_1}\in\mathbb{N}^+, \exists~ M_2=M_2(a^1_{M_1})\in\mathbb{N}^+, s.t. \\
& \qquad \cdots\cdots ( \text{recursively})\\
& \quad\qquad \forall~ a^{m-1}_0<a^{m-1}_1<\cdots<a^{m-1}_{M_{m-1}}\in\mathbb{N}^+, \exists~ M_m=M_m(a^{m-1}_{M_{m-1}})\in\mathbb{N}^+, s.t. \\
&\qquad\qquad\Ord A(X,d)^{\big\{ \left. a^i_{p_i} \right|~0\le i\le m-1, 0\le p_i\le M_i \big\}}\le M_m.
\end{align*}
Since $A(X,d)$ is inclusive (i.e., for every $\sigma, \tau \in \mathrm{Fin},\mathbb{N}$ with $\tau \subseteq \sigma$, if $\sigma \in A(X,d)$ then $\tau \in A(X,d)$), we have
\begin{align*}
&\Ord A(X,d)^{\big\{ a^i_{p_i}~\mid ~0\le i\le m-1, 0\le p_i\le M_i \big\}}\le M_m\\
\iff &\forall~ a^{m}_0<a^{m}_1<\cdots<a^{m}_{M_{m}}\in\mathbb{N}^+, \\
&\big\{ \left. a^i_{p_i} ~\right|~0\le i\le m-1, 0\le p_i\le M_i \big\}\cup \big\{a^{m}_0, a^{m}_1, \cdots, a^{m}_{M_{m}} \big\}\notin A(X,d). 
\end{align*}
We now proceed to prove the equivalence.

``$(1)\Longrightarrow (2)$'': Given $k_0, k_1, \cdots, k_m \in \mathbb{N}^+$ with $k_0 < k_1 < \cdots < k_m$, define $a^0_p := k_0 + p$ for all $0 \le p \le n$. By the equivalence discussed above, there exists an integer $M_1$ depending on $k_0 + n$, and we set $E_1 := M_1(k_0 + n)$. 
Then, for each $1 \le i \le m-1$, recursively define $a^i_{p_i} := k_i + p_i$ for all $0 \le p_i \le M_i$. During this recursive process, we again apply the same equivalence to obtain the existence of $M_{i+1}$ depending on $a^i_{M_i}$, and we define $E_{i+1} := M_{i+1}(a^i_{M_i})$. 
Finally, define $a^m_{p_m} := k_m + p_m$ for all $0 \le p_m \le M_m$. Then, as discussed above, it follows that there exist families satisfying the required conditions.

``$(2)\Longrightarrow (1)$'': For any increasing sequence $a^0_0 < a^0_1 < \cdots < a^0_n \in \mathbb{N}^+$, define $k_0 := a^0_n$ and let $M_1(a^0_n) := E_1$. Next for any $a^1_0 < a^1_1 < \cdots < a^1_{M_1} \in \mathbb{N}^+$, define $k_1 := \max\{k_0, a^1_{M_1}\} + 1$, and let $M_2(a^1_{M_1}) := E_2$. Proceeding recursively, for each $1 \le i \le m-1$, define 
\[
k_i := \max\{k_0, k_1, \cdots, k_{i-1}, a^i_{M_i}\} + 1,\text{ and let } M_{i+1}(a^i_{M_i}) := E_{i+1}.
\]
Finally, consider $k_m := \max\{k_0, k_1, \cdots, k_{m-1}, a^m_{M_m}\} + 1$. Then by assumption (2), we conclude that $\trasdim(X) \le \omega {m} + n$.
\end{proof}

We then strengthen the condition $\trasdim(X) \leq \omega {m}+n$ by additionally requiring a form of uniform disjointness under translation.

\subsection{A Technical Tool: Disjointness of the Tiling through Uniform Translation}

\begin{definition}\label{def:(m,n)-DRUT}
Let $(X,d)$ be a metric space, let $h\colon \mathbb{N}^+ \to \mathbb{N}^+$ be a map, and let $m, n \in \mathbb{N}$.
We say that $X$ has the property of \emph{disjointness of the $(m,n)$-tiling through uniform translation associated with $h$} if there exist integers $F_1, F_2, \cdots, F_m \in \mathbb{N}$, depending only on $X$ and independent of the choice of $h$, such that for every sequence $k_0, k_1, \cdots, k_m \in \mathbb{N}^+$ with $k_0 < k_1 < \cdots < k_m$, and for each vector $\vec{l}=(l_1,\cdots,l_m)\in \prod_{r=1}^m \{1,2,\cdots,h(k_r)\}$, there exist uniformly bounded families
\[
\mathcal{C}_{0}^{\vec{l}},\, \mathcal{C}_{1}^{\vec{l}},\, \cdots,\, \mathcal{C}_{n}^{\vec{l}},\quad 
\mathcal{D}_{i,0}^{\vec{l}},\, \mathcal{D}_{i,1}^{\vec{l}},\, \cdots,\, \mathcal{D}_{i,F_i}^{\vec{l}}\quad (1\le i\le m)
\]
such that for every fixed $\vec{l}$, the following holds: 
\begin{itemize}
    \item Each $\mathcal{C}_j^{\vec{l}}$ is $k_0$-disjoint for all $0\le j \le n$.
    \item Each $\mathcal{D}_{i,s_i}^{\vec{l}}$ is $k_i$-disjoint for all $1\le i\le m$ and $0\le s_i \le F_i$.
    \item The family $\big(\bigcup_{j=0}^{n}\mathcal{C}_j^{\vec{l}}\big)\cup\big(\bigcup_{i=1}^m\bigcup_{s_i=0}^{F_i}\mathcal{D}_{i,s_i}^{\vec{l}}\big)$ covers $X$.
\end{itemize}

Moreover, for each $1\le i\le m$, $0\le s_i\le F_i$, and each fixed $(l_1,\cdots, l_{i-1},l_{i+1},\cdots,l_m)$, the following holds: 
\begin{itemize}
    \item For all $1\le t_i< t_i' \le h(k_i)$, the sets $\bigcup \mathcal{D}_{i,s_i}^{l_1,\cdots, l_{i-1}, t_i, l_{i+1} \cdots,l_m}$ and $\bigcup \mathcal{D}_{i,s_i}^{l_1,\cdots, l_{i-1}, t_i', l_{i+1} \cdots,l_m}$ are non-empty and disjoint. 
    \item The family $\bigcup_{t_i=1}^{h(k_i)}\mathcal{D}_{i,s_i}^{l_1,\cdots, l_{i-1}, t_i, l_{i+1} \cdots,l_m}$ is $k_i$-disjoint.
\end{itemize}

We say that $X$ has the property of \emph{disjointness of the $(m,n)$-tiling through uniform translation} (abbreviated as $(m,n)$-DTUT) if the above condition holds for all maps $h\colon \mathbb{N}^+\to\mathbb{N}^+$. 

Furthermore, we say that $X$ has the property of \emph{disjointness of the tiling through uniform translation} (abbreviated as DTUT) if there exist $m, n \in \mathbb{N}$ such that $X$ has the $(m,n)$-DTUT property.
\end{definition}
Clearly, if $X$ has the $(m,n)$-DTUT property, then $\trasdim(X) \le \omega {m} + n$.

Recall that a group $\Gamma$ is called an \emph{FC-group} if for all $\gamma\in\Gamma$, the conjugacy class of $\gamma$ is finite. 
\begin{lemma}\label{lem:fgFChasFAD}
Finitely generated FC-groups have finite asymptotic dimension.
\end{lemma}
\begin{proof}
Suppose $\Gamma$ is an FC-group with a symmetric generating set $S=\{s_1,\cdots,s_m\}$. Denote the center of $\Gamma$ by $Z(\Gamma)$, then we have $Z(\Gamma)=\bigcap_{s_i\in S}Z(s_i)$, where $Z(s_i)=\{g\in \Gamma \mid gs_i=s_ig\}$ is the centralizer of $s_i$ for all $s_i\in S$. 
Since the index of $Z(s_i)$ in $\Gamma$ equals the number of its conjugacy classes, which is finite, the index of $Z(\Gamma)$ is also finite. 

Note that in a finitely generated group, the subgroup with finite index is also finitely generated. Therefore $Z(\Gamma)$ is a finitely generated abelian group, moreover, $Z(\Gamma)$ has finite asymptotic dimension, and we obtain that $\Gamma$ has finite asymptotic dimension through the permanence property. 
\end{proof}

Throughout this paper, let $H$ denote an FC-group equipped with a proper left-invariant metric, and suppose that $\asdim(H) = h < \infty$. For each $k>0$ and $0\le p\le h$, let $^{k}\mathcal{H}_p:=\{^{k}H_{pq} \mid q\in\mathbb{N}\}$ be a $D(k)$-bounded family of subsets of $H$, such that the family $\bigcup_{p=0}^h {^k\mathcal{H}_p}$ covers $H$, and each ${^k\mathcal{H}_p}$ is $k$-disjoint. 

Given $W\in \WR(G, H\curvearrowright I)$, where $G$ is a group with finite asymptotic dimension and satisfies the $(m,0)$-DTUT property, and $H$ is an FC-group with finite asymptotic dimension, we aim to decompose $W$ in such a way that $\trasdim(W) \leq \omega {m} + \asdim (H)$. Note that for any $k\in\mathbb{N}^+$, we have 
\[
W= \bigg( \bigoplus_{I}G, H \bigg)=\bigcup_{p=0}^h\bigcup_{q\in\mathbb{N}}\bigg( \bigoplus_{I}G, {^{k}H_{pq}} \bigg).
\] 
Then, for any sequence $k_0, k_1, \cdots, k_m \in \mathbb{N}^+$ with $k_0 < k_1 < \cdots < k_m$, our general strategy is as follows: it suffices to first decompose each subset of the form $\big( \bigoplus_{I}G, {^{k_0}H_{pq}} \big)$ for all $0 \leq p \leq h$ and $q \in \mathbb{N}$ in a sufficiently fine manner.
After that, by repeatedly applying the saturated union technique, we obtain decompositions for subsets of the form $\big( \bigoplus_{I}G, {^{k_m}H_{pq}} \big)$. This ultimately enables us to construct a covering of $W$ satisfying the disjointness conditions required for transfinite asymptotic dimension, thus verifying asymptotic property C. 

\subsection{A Technical Lemma}
We have the following technical lemma for decomposing each $\big( \bigoplus_{I}G, {^{k_0}H_{pq}} \big)$.  

\begin{lemma}\label{lemma: decomp for the infinite sum (pre for wr prod)}
Let $G$ be a group with finite asymptotic dimension that satisfies the $(m,0)$-DTUT property for some $m\in\mathbb{N}^+$, and let $H$ be an FC-group with finite asymptotic dimension. 
Let $W\in \WR(G, H\curvearrowright I)$ with finitely many orbits $\{I_\alpha\}_{\alpha\in \orb(I)}$. Suppose $T$ is a symmetric transversal associated with ${G}^{(I)}\lhd W$ such that the transversal map $\eta: H\longrightarrow T$ is unital and bornologous. 

For every sequence $k_0,k_1,\cdots,k_m\in\mathbb{N}^+$ with $k_0< k_1<\cdots<k_m$, for every $0\le p\le h$ and $q\in\mathbb{N}$, denote $\displaystyle X_{pq}:=(\bigoplus_{I}G, {^{k_0}H_{pq}})$, then, there exists $B(k_m)>0$ and $B(k_m)$-bounded families of subsets of $W$
    \begin{align*}
    \mathcal{U}_{0},\ \mathcal{U}_{i,j_i,s,t_i} 
    \text{ with } &
    1\le i\le m,~ 0\le j_i\le F_i,~s\in T(k_0), \text{ and }\\
    &1\le t_i\le (g+1)^{|T(k_{i-1})|},
    \end{align*}
such that $\mathcal{U}_0$ is $k_0$-disjoint, each $\mathcal{U}_{i,j_i,s,t_i}$ is $k_i$-disjoint, and $\big(\bigcup_{i,j_i,s,t_i}\mathcal{U}_{i,j_i,s,t_i}\big)\cup\mathcal{U}_0$ covers $X_{pq}$. 
Here, $F_i$ denotes the integer appearing in the definition of the $(m,0)$-DTUT property; $g$ denotes the asymptotic dimension of $G$, and for each $1 \le j \le m$, we denote 
\[
T(k_j) := \bigcup_{\alpha \in \orb(I)} \big({^{k_0}H_{pq}^{-1}} B_H(1_H, k_j + P(k_j)) \big){.} i_{\alpha,0} \subseteq I,
\]
where $P:\mathbb{R}\to\mathbb{R}$ is a control function associated with the bornologous map $\eta$.
\end{lemma}

Note that for all $1 \le i \le m$, the number of families $\mathcal{U}_{i,j_i,s,t_i}$ obtained by ranging over $j_i$, $s$, and $t_i$ depends on the pair $\{k_0, k_{i-1}\}$.

\begin{proof}
Since $\asdim(G)=g<\infty$, for all $u>0$, and $0\le v\le g$, let $^{u}\mathcal{G}_v:=\{^u{G}_{vr} \mid r\in\mathbb{N}\}$ be a $V(u)$-bounded family of subsets of $G$, such that the family $\bigcup_{v=0}^g {^u\mathcal{G}_v}$ covers $G$, and each ${^u\mathcal{G}_v}$ is $u$-disjoint. 
Let $P:\mathbb{R}\to\mathbb{R}$ be a control function associated with the map $\eta: H\longrightarrow T$ that is unital and bornologous, as described in Remark~\ref{remark:bornologous}.

For every sequence $k_0, k_1, \cdots, k_m \in \mathbb{N}^+$ with $k_0 < k_1 < \cdots < k_m$, we may assume without loss of generality that $P(k_0) > |\orb(I)|$. Define
\begin{align*}
Q(k_j):=k_j+\sum_{i=1}^{k_j} P(i)+ \sum_{j=1}^{P(D(k_j))}2|B_H(1_H,j)|\times \big(P(D(k_j))\big)^{k_j}\times |\orb(I)| 
\end{align*}
for all $0\le j\le m$. 
First, note that for any $w=((a_i)_{i\in I},h)\in W$, if $l_W(w)\le k_j$ for some $0\le j\le m$, then
\[
\supp((a_i)_{i\in I})\subset \bigcup_{\alpha\in \orb(I)}\big(B_H(1_H,P(k_j))\big){.} i_{\alpha,0}, 
\]
and for any $j\in I$, 
\[
l_G(a_j)\le k_j+\sum_{i=1}^{k_j} P(i)<Q(k_j)\le Q(k_m).
\]

Next, note that for any $\alpha,\beta\in {^{k_0}H_{pq}}$, we have $d_H(\alpha, \beta)\le D(k_0)$, hence
\begin{align*}       &l_W\bigg(\big(\eta(\alpha^{-1})\eta(\beta)\eta(\alpha^{-1}\beta)^{-1}, \alpha^{-1}\beta\big) \bigg)=l_W\big( (e,\alpha^{-1})\cdot (e,\beta) \big)\\
\le &\sum_{j=1}^{P(D(k_0))}2|B_H(1_H,j)|\times P(D(k_0))^2\times |\orb(I)|\\
\le &Q(k_0). 
\end{align*}

For each $1 \le i \le m$, define
\begin{align*}
T(k_i,k_{i-1}):=T(k_i)\setminus T(k_{i-1})\subset I.
\end{align*} 

Fix a map $h \colon \mathbb{N}^+ \to \mathbb{N}^+$ such that $h(Q(k_i)) = (g+1)^{|T(k_i)|}$ for all $1 \le i \le m$. 
Since $G$ has the $(m,0)$-DTUT property (associated with $h$), there exist integers $F_1, F_2, \cdots, F_m \in \mathbb{N}$, such that for the sequence $Q(k_0) < Q(k_1) < \cdots < Q(k_m)$, and for each vector $\vec{l}=(l_1,\cdots,l_m) \in \prod_{r=1}^m \{1,2,\cdots,(g+1)^{|T(k_r, k_{r-1})|}\}$, there exist $R(Q(k_m))$-bounded families
\[
\mathcal{C}_{0}^{\vec{l}},\quad 
\mathcal{D}_{i,0}^{\vec{l}},\, \mathcal{D}_{i,1}^{\vec{l}},\, \cdots,\, \mathcal{D}_{i,F_i}^{\vec{l}} \quad (1 \le i \le m)
\]
such that these families satisfy the corresponding disjointness and covering conditions in Definition~\ref{def:(m,n)-DRUT}.

Let $\phi_0 \colon \{1,2,\cdots,(g+1)^{|T(k_0)|}\} \to \{0,1,\cdots,g\}^{T(k_0)}$ be a fixed bijection. For each $1 \le i \le m$, let 
\[
\phi_i \colon \{1,2,\cdots,(g+1)^{|T(k_i,k_{i-1})|}\} \longrightarrow \{0,1,\cdots,g\}^{T(k_i,k_{i-1})}
\]
be a fixed bijection. Then, for all $l_i \in \{1,2,\cdots,(g+1)^{|T(k_i,k_{i-1})|}\}$, define
\[
\mathcal{W}_{l_i} := \bigg\{ \prod_{w \in T(k_i,k_{i-1})} V_w \mid V_w \in {^{Q(k_m)}\mathcal{G}_{\phi_i(l_i)_w}} \bigg\}.
\]
Note that the family $\bigcup_{l_i=1}^{(g+1)^{|T(k_i,k_{i-1})|}} \mathcal{W}_{l_i}$ is pairwise disjoint, and covers $\bigoplus_{T(k_i,k_{i-1})}G$. 
Now we define 
\begin{align*}
\mathcal{U}_{0}:=\bigg\{\big(&\prod_{i\in T(k_0)}C_i\times W_1\times \cdots\times W_m \times \{\widetilde{x}\}, {^{k_0}H_{pq}}\big)~\bigg|~  \\ 
&\widetilde{x}\in \bigoplus_{I}G \text{ and }\supp(\widetilde{x})\cap T(k_m) =\varnothing, \\
&W_j \subset \prod_{T(k_j,k_{j-1})}G \quad (1\le j\le m), \\
&\prod_{i\in T(k)}C_i\times W_1\times \cdots\times W_m \in\bigcup_{\substack{1\le n\le m \\ l_n=1}}^{(g+1)^{|T(k_n,k_{n-1})|}} \prod_{T(k_0)}\mathcal{C}_0^{l_1,\cdots,l_m}\times \mathcal{W}_{l_1}\times\cdots\times\mathcal{W}_{l_m}
\bigg\}.
\end{align*}

Note that for all $t_i \in \{1, 2, \cdots, (g+1)^{|T(k_{i-1})|}\}$ ($1 \le i \le m$), we have $|T(k_{i-1})| = \sum_{r=0}^{i-1} |T(k_r, k_{r-1})|$, where we interpret $T(k_0, k_{-1})$ as $T(k_0)$. Therefore, we may associate each $t_i$ uniquely with a vector $\vec{t_i} = (t_{i,0}, t_{i,1}, \cdots, t_{i,i-1})$, where for all $0 \le r \le i-1$, we have
\[
t_{i,r} \in \{1, 2, \cdots, (g+1)^{|T(k_r, k_{r-1})|}\}.
\]
This correspondence between $t_i$ and $\vec{t_i}$ is one-to-one. Then for all $1\le i\le m$, $0\le j_i\le F_i$, $s\in T(k_0)$, and $t_i\in\{1,2,\cdots,(g+1)^{|T(k_{i-1})|}\}$, we define 
\begin{align*}
\mathcal{U}_{i,j_i,s,t_i}:=
\bigg\{\big(&D_{s}\times\prod_{u\in T(k_{0})\setminus\{s\}}V_u\times Y_1\times\cdots Y_{i-1}\times W_i\times\cdots W_m\times \{\widetilde{x}\}, {^{k_0}H_{pq}}\big)~\bigg|~  \\
&\widetilde{x}\in \bigoplus_{I}G \text{ and }\supp(\widetilde{x})\cap T(k_m) =\varnothing, \\ 
&Y_j \subset \prod_{T(k_j,k_{j-1})}G \quad (1\le j\le i-1), \quad W_n \subset \prod_{T(k_j,k_{j-1})}G \quad (i\le n\le m), \\
&V_u\in {^{Q(k_m)}\mathcal{G}_{\phi_0(t_{i,0})_{u}}} \quad (u\in T(k_0)\setminus\{s\}),\\
&D_{s}\times Y_1\times\cdots Y_{i-1}\times W_i\times\cdots W_m\in\\
&\bigcup_{\substack{i\le n\le m \\ l_n=1}}^{(g+1)^{|T(k_n,k_{n-1})|}} \mathcal{D}_{i,j_i}^{t_{i,1},\cdots,t_{i,i-1},l_i,\cdots,l_m}\times \mathcal{W}_{t_{i,1}}\times \cdots \times\mathcal{W}_{t_{i,i-1}}\times \mathcal{W}_{l_i}\times\cdots\times\mathcal{W}_{l_m} 
\bigg\}.
\end{align*}

Now we check that these families of subsets of $X_{pq}$ satisfy the required conditions with the following choice of $B(k_m)$: 
\begin{align*}
B(k_m)=&|\orb(I)|\times [R(Q(k_m))+V(Q(k_m))]\\
     &~ \times|B_H(1_H,k_m+P(k_m))|\times|B_H(1_H,D(k_m))|\\
     &~ +2|\orb(I)|\times \big(k_m+P(k_m)+M_H(k_m)\big)^2\big|B_H\big(1_H,k_m+P(k_m)+M_H(k_m)\big)\big|\\
     &~ +Q(k_m),
\end{align*}
where 
\[
M_H(k_m):=\sup_{h\in H}\{l_H(hgh^{-1}) \mid g\in B_H(1_H,D(k_m))\}
\]
is finite, since $H$ is an FC-group. 
We provide the detailed verification steps as follows.

\noindent\textbf{Step 1. }Prove that $\mathcal{U}_{0}$ is $k_0$-disjoint and $B(k_m)$-bounded.
\begin{enumerate}
  \item Let 
        \begin{align*}
        &U=\big(\prod_{i\in T(k_0)}C_i\times W_1\times \cdots\times W_m \times \{\widetilde{x}\}, {^{k_0}H_{pq}}\big),\\
        &U'=\big(\prod_{i\in T(k_0)}C_i'\times W_1'\times \cdots\times W_m' \times \{\widetilde{x}'\}, {^{k_0}H_{pq}}\big),
        \end{align*}
        and $U\neq U'$.
        Then for all $x=(f,\alpha)\in U$, and $y=(g,\beta)\in U'$, we distinguish two cases as follows:
        \begin{itemize}
          \item If $( W_1\times \cdots\times W_m \times \{\widetilde{x}\})\neq ( W_1'\times \cdots\times W_m' \times \{\widetilde{x}'\})$, then there exists $j\notin T(k_0)$ such that $f(j)\neq g(j)$, i.e., $(f^{-1}g)(\alpha^{-1}\cdot \alpha{.} j)\neq 1_{G}$. 
          
          Since $\alpha{.} j \notin \bigcup_{\gamma\in \orb(I)}\big(B_H(1_H,k_0+P(k_0))\big){.} i_{\gamma,0}$, we have $l_W(x^{-1}y)\ge k_0$ and $d_W(U,U')\ge k_0$. 
          
          \item If $(W_1 \times \cdots \times W_m \times \{\widetilde{x}\}) = (W_1' \times \cdots \times W_m' \times \{\widetilde{x}'\})$, suppose both $W_1 \times \cdots \times W_m$ and $W_1' \times \cdots \times W_m'$ belong to $\mathcal{W}_{l_1} \times \cdots \times \mathcal{W}_{l_m}$ for the same $(l_1, \cdots, l_m)$. Then there exists $j \in T(k_0)$, with $C_j, (C_j)' \in \mathcal{C}_0^{l_1,\cdots,l_m}$ and $C_j \neq (C_j)'$. Therefore, we have $d_G(f(j), g(j)) \ge Q(k_0)$, which implies $d_W(x, y) \ge k_0$ and hence $d_W(U, U') \ge k_0$. 
        \end{itemize}
  \item Let $x=(f,\alpha), z=(u,\xi)\in U$, we consider the following cases: 
        \begin{itemize}
          \item If $f=u$, then $d_W(x,z)=l_W\big((e,\alpha^{-1})\cdot (e,\xi)\big)\le Q(k_m)$. 
          \item If $f\neq u$, then 
                \begin{align*}
                \sum_{i\in I}d_{G}(f(i),u(i))&\le |T(k_0)|\times R(Q(k_m))+\sum_{j=1}^m|T(k_j,k_{j-1})|\times V(Q(k_m))\\
                &\le|\orb(I)|\times|B_H(1_H,k_m+P(k_m))|\times|B_H(1_H,D(k_m))|\\
                &\quad\times [R(Q(k_m))+V(Q(k_m))].
                \end{align*}
                Note that ${^{k_0}H_{pq}}{^{k_0}H_{pq}^{-1}}\subset B_H(1_H,\text{diam }{^{k_0}H_{pq}^{-1}})$, therefore 
                \begin{align*}
                &\supp\bigg(\big((f^{-1}u)_{\alpha^{-1}{.} i}\big)_{i\in I}\bigg) \subset  {^{k_0}H_{pq}}\cdot T(k_m)\\
                \subset & \bigcup_{\alpha\in \orb(I)}B_H\big(1_H,k_m+P(k_m)+M_H(k_m)\big){.} i_{\alpha,0}.
                \end{align*}
                It follows that 
                \begin{align*}
                &d_W(x,z)=l_W\bigg( \big((f^{-1}u)_{\alpha^{-1}{.} i}\big)_{i\in I}\cdot \eta(\alpha^{-1})\eta(\xi)\eta(\alpha^{-1}\xi ), ~\alpha^{-1}\xi  \bigg)\\ 
                \le & l_W\big( \big((f^{-1}u)_{\alpha^{-1}{.} i}\big)_{i\in I}, 1_H \big)+l_W\big( (e,\alpha^{-1})\cdot (e,\beta) \big)\\ 
                \le & 2|\orb(I)|\times \big(k_m+P(k_m)+M_H(k_m)\big)^2\times \big|B_H\big(1_H,k_m+P(k_m)+M_H(k_m)\big)\big|\\ 
               &+\sum_{i\in I}d_{\mathbb{F}_2}(f(i),u(i))+Q(k_m)\\
               \le & B(k_m) .
                \end{align*}
        \end{itemize}
\end{enumerate}

\noindent\textbf{Step 2. }Prove that for all $1\le i\le m$, $0\le j_i\le F_i$, $s\in T(k_0)$, and $1\le t_i\le (g+1)^{|T(k_{i-1})|}$, the family $\mathcal{U}_{i,j_i,s,t_i}$ is $k_i$-disjoint and $B(k_m)$-bounded.
\begin{enumerate}
  \item Let
        \begin{align*}
        &U=\big(D_{s}\times\prod_{u\in T(k_{0})\setminus\{s\}}V_u\times Y_1\times\cdots \times  Y_{i-1}\times W_i\times\cdots \times W_m\times \{\widetilde{x}\}, {^{k_0}H_{pq}}\big),\\
        &U'=\big(D_{s}'\times\prod_{u\in T(k_{0})\setminus\{s\}}V_u'\times Y_1'\times\cdots \times Y_{i-1}'\times W_i'\times\cdots \times W_m'\times \{\widetilde{x}'\}, {^{k_0}H_{pq}}\big).
        \end{align*}
        and $U\neq U'$. 
        Then for all $x=(f,\alpha)\in U$ and $y=(g,\beta)\in U'$, we consider the following cases:
        \begin{itemize}
           \item If $(W_{i+1} \times \cdots \times W_m \times \{\widetilde{x}\}) \neq (W_{i+1}' \times \cdots \times W_m' \times \{\widetilde{x}'\})$ (note that if $i = m$, this reduces to $\{\widetilde{x}\} \neq \{\widetilde{x}'\}$), then there exists some $j \notin T(k_i)$ such that $f(j) \neq g(j)$, i.e., $(f^{-1}g)(\alpha^{-1} \cdot \alpha{.} j) \neq 1_G$. Since $\alpha{.} j \notin \bigcup_{\gamma \in \orb(I)} B_H(1_H, k_i + P(k_i)){.}i_{\gamma,0}$, we conclude that $l_W(x^{-1} y) \ge k_i$, and thus $d_W(U, U') \ge k_i$.

           \item If $(W_{i+1} \times \cdots \times W_m \times \{\widetilde{x}\}) = (W_{i+1}' \times \cdots \times W_m' \times \{\widetilde{x}'\})$, then suppose that both $W_{i+1} \times \cdots \times W_m$ and $W_{i+1}' \times \cdots \times W_m'$ belong to $\mathcal{W}_{l_{i+1}} \times \cdots \times \mathcal{W}_{l_m}$ for the same $(l_{i+1}, \cdots, l_m)$. In this case, we have
           \begin{align*}
            D_s, D_s' \in \bigcup_{l_i = 1}^{(g + 1)^{|T(k_i, k_{i-1})|}} \mathcal{D}_{i, j_i}^{t_{i,1}, \cdots, t_{i,i-1}, l_i, l_{i+1} \cdots, l_m}.
           \end{align*}
           Then we distinguish the following subcases:
            \begin{itemize}
                \item If $D_s \neq D_s'$, since 
                $\bigcup_{l_i = 1}^{(g + 1)^{|T(k_i, k_{i-1})|}} \mathcal{D}_{i, j_i}^{t_{i,1}, \cdots, t_{i,i-1}, l_i, l_{i+1}, \cdots, l_m}$
                is $Q(k_i)$-disjoint, we have $d_W(x, y) \ge k_i$.
            
                \item If $D_s = D_s'$, suppose both $D_s$ and $D_s'$ belong to 
                $\mathcal{D}_{i, j_i}^{t_{i,1}, \cdots, t_{i,i-1}, l_i, l_{i+1}, \cdots, l_m}$ for the same $(l_i, l_{i+1}, \cdots, l_m)$. In this case, we have $W_i, W_i' \in \mathcal{W}_{l_i}$, and
                \[
                \prod_{u \in T(k_0) \setminus \{s\}} V_u \times Y_1\times \cdots \times Y_{i-1}\times W_i \neq \prod_{u \in T(k_0) \setminus \{s\}} V_u' \times Y_1'\times \cdots \times Y_{i-1}'\times W_i'.
                \]
                It follows that at least one of the following two situations holds:
                \begin{itemize}
                    \item There exists $u \in T(k_0) \setminus \{s\}$ such that $V_u \neq V_u'$, which implies that $d_G(f(u), g(u)) \ge Q(k_m)$, and hence $d_W(x, y) \ge k_m$.
                    \item In the case when $Y_1 \times \cdots \times Y_{i-1} \times W_i \ne Y_1' \times \cdots \times Y_{i-1}' \times W_i'$, since $Y_p, Y_p' \in \mathcal{W}_{t_{i,p}}$ for all $1 \le p \le i-1$ and $W_i, W_i' \in \mathcal{W}_{l_i}$, there must exist $j \in T(k_i) \setminus T(k_0)$ such that $d_G(f(j), g(j)) \ge Q(k_m)$, and thus $d_W(x, y) \ge k_m$.
                \end{itemize}
            \end{itemize}
        \end{itemize}
        Therefore, we have $d_W(U,U')\ge k_i$.
  \item Let $x=(f,\alpha), z=(u,\xi)\in U$, then by the same argument as in \textbf{Step 1}, we obtain $d_W(x,z) \le B(k_m)$.
\end{enumerate}

\noindent\textbf{Step 3. }Prove that $\big(\bigcup_{i,j_i,s,t_i}\mathcal{U}_{i,j_i,s,t_i}\big)\cup\mathcal{U}_0$ covers $X_{pq}$, where $1\le i\le m$, $0\le j_i\le F_i$, $s\in T(k_0)$, and $1\le t_i \le (g+1)^{|T(k_{i-1})|}$.

Let $x = (f, \alpha) \in X_{pq} \setminus \bigcup \mathcal{U}_0$. Suppose that the tuple $\big(f(r)\big)_{r \in T(k_m)\setminus T(k_0)}$ belongs to $\mathcal{W}_{l_1}\times\cdots\times\mathcal{W}_{l_m}$ for some $\vec{l}=(l_1,\cdots,l_m) \in \prod_{r=1}^m \{1,2,\cdots,(g+1)^{|T(k_r, k_{r-1})|}\}$. 
Since $x \notin \bigcup \mathcal{U}_0$, it follows that 
\[
(f(j))_{j \in T(k_0)} \notin \bigcup \left\{ \prod_{r \in T(k_0)} C_r \mid C_r \in \mathcal{C}_0^{\vec{l}} \right\}.
\]
Therefore, there exists $s \in T(k_0)$ such that $f(s) \notin \bigcup \mathcal{C}_0^{\vec{l}}$.
Since $\mathcal{C}_0^{\vec{l}}\cup\big(\bigcup_{i=1}^m\bigcup_{j_i=0}^{F_i}\mathcal{D}_{i,j_i}^{\vec{l}}\big)$ covers $G$, then there exists $1\le i\le m$, $0\le j_i\le F_i$, and $D_s \in \mathcal{D}_{i,j_i}^{\vec{l}}$ such that $f(s) \in D_s$. Consequently, we can choose $1\le t_i \le  2^{|T(k_{i-1})|}$ such that $\vec{t_i} = (t_{i,0}, t_{i,1}, \cdots, t_{i,i-1})$ satisfies that $t_{i,v}=l_v$ for all $1\le v\le i-1$, and 
\[
(f(j))_{j \in T(k_0)} \in \bigcup \left\{ \prod_{u \in T(k_0)} V_u \mid V_u \in {^{Q(k_m)}}\mathcal{G}_{\phi_0(t_{i,0})_u} \right\}.
\]
It follows that $x \in \bigcup \mathcal{U}_{i,j_i,s,t_i}$. 
\end{proof}

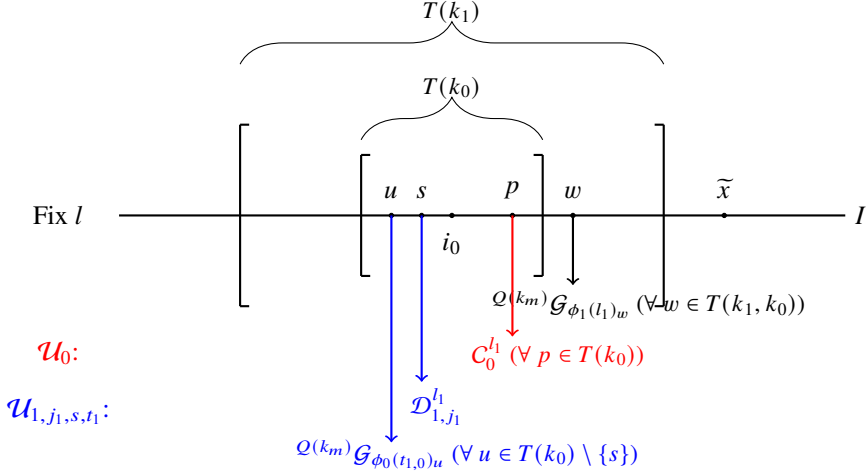
\begin{figure}[h]
\centering
\begin{tikzpicture}[scale=0.4]

\draw[thick] (-11,0) -- (13,0);

\draw[decorate,decoration={brace,amplitude=16pt}] (-3,2.5) -- (3,2.5) node[midway,yshift=12pt,above] {\small$T(k_0)$};
\draw[decorate,decoration={brace,amplitude=16pt}] (-7,5) -- (7,5) node[midway,yshift=12pt,above] {\small$T(k_1)$} ;

\node at (0, 0) [below=3pt] {$i_0$};
\draw[fill=black] (0, 0) circle (2pt); 

\node at (-1, 0) [above=3pt] {$s$};
\draw[fill=black] (-1, 0) circle (2pt); 

\node at (-2, 0) [above=3pt] {$u$};
\draw[fill=black] (-2, 0) circle (2pt); 

\node at (2, 0) [above=3pt] {$p$};
\draw[fill=black] (2, 0) circle (2pt); 

\node at (4, 0) [above=3pt] {$w$};
\draw[fill=black] (4, 0) circle (2pt);

\node at (9, 0) [above=3pt] {$\widetilde{x}$};
\draw[fill=black] (9, 0) circle (2pt);

\node at (-13, 0) [] {$\text{Fix }l$};
\node at (13.5, 0) [] {$I$};
\node[text=red] at (-13, -5.5) [above=3pt] {$\mathcal{U}_{0}$: };
\node[text=blue] at (-13, -7.5) [above=3pt] {$\mathcal{U}_{1,j_1,s,t_1}$: };

\node[text=blue] at (0.5, -9) [above=3pt] {\small ${^{Q(k_m)}\mathcal{G}_{\phi_0(t_{1,0})_u}}$ ($\forall~u\in T(k_0)\setminus \{s\}$)};
\node[text=blue] at (-0.5, -7.5) [above=3pt] {\small $\mathcal{D}_{1,j_1}^{l_1}$};
\node[text=red] at (3.5, -5.7) [above=3pt] {\small $\mathcal{C}_0^{l_1}$ ($\forall~p\in T(k_0)$)};
\node at (6.5, -4) [above=3pt] {\small ${^{Q(k_m)}\mathcal{G}_{\phi_1(l_1)_w}}$ ($\forall~w\in T(k_1,k_0)$)};

\draw[->, thick, blue] (-2, 0) -- (-2, -7.5) node[midway, left] {};
\draw[->, thick, blue] (-1, 0) -- (-1, -5.5) node[midway, left] {};
\draw[->, thick, red] (2, 0) -- (2, -4) node[midway, left] {};
\draw[->, thick] (4, 0) -- (4, -2.3) node[midway, left] {};

\draw[thick] (-7,3) -- (-7,-3);
\draw[thick] (-7,3) -- (-6.7,3);
\draw[thick] (-7,-3) -- (-6.7,-3);
\draw[thick] (7,3) -- (7,-3);
\draw[thick] (7,3) -- (6.7,3);
\draw[thick] (7,-3) -- (6.7,-3);

\draw[thick] (-3,2) -- (-3,-2);
\draw[thick] (-3,2) -- (-2.7,2);
\draw[thick] (-3,-2) -- (-2.7,-2);
\draw[thick] (3,2) -- (3,-2);
\draw[thick] (3,2) -- (2.7,2);
\draw[thick] (3,-2) -- (2.7,-2);
\end{tikzpicture}
\caption{Ideas for the decomposition method for $\big( \bigoplus_{I}G, {^{k_0}H_{pq}} \big)$.}
\label{fig:ideas_for_the_constructions_in_lemmas}
\end{figure}

Figure~\ref{fig:ideas_for_the_constructions_in_lemmas} shows the decomposition of $\big( \bigoplus_{I}G, {^{k_0}H_{pq}} \big)$ in the case where the group $G$ satisfies the $(1,0)$-DTUT property. In the figure, the indices range over $0\le j_1 \le  F_1$, $s \in T(k_0)$, and $1\le t_1 \le (g+1)^{|T(k_0)|}$, where the associated vector is given by $\vec{t_1} = (t_{1,0})$. 

\subsection{The Main Result}
We aim to use the decomposition given in \cref{lemma: decomp for the infinite sum (pre for wr prod)}, and by repeatedly applying the saturated union technique, to construct a covering of $W$ that satisfies the disjointness conditions required for the transfinite asymptotic dimension.

\begin{theorem}\label{thm: GwrlikeH has apc}
Let $G$ be a group with finite asymptotic dimension that satisfies the $(m,0)$-DTUT property for some $m\in\mathbb{N}$.
Let $H$ be an FC-group with finite asymptotic dimension.
Let $W\in \WR(G, H\curvearrowright I)$ with finitely many orbits $\{I_\alpha\}_{\alpha\in \orb(I)}$. 
Suppose $T$ is a symmetric transversal associated with ${G}^{(I)}\lhd W$ such that the transversal map $\eta: H\longrightarrow T$ is unital and bornologous. 
Then  $\trasdim(W) \le \omega \,{m} + \asdim(H)$, and hence $W$ has asymptotic property~C. 
\end{theorem}

\begin{proof}
As stated before, we assume $\asdim(H) = h < \infty$ throughout. 
Consider the case when $m = 0$. Then $\big( \bigoplus_{I} G, d_W \big)$ has finite asymptotic dimension, and hence $W$ has asymptotic property C by \cref{thm: closed properties for APC}. Therefore, we may assume that $G$ has the $(m,0)$-DTUT property with $m > 0$.

To simplify notation, we denote $|B_H(1_H, D(k_i))|$ by $D_H(k_i)$ for each $0 \le i \le m$.

For every sequence $k_0,k_1,\cdots,k_m\in\mathbb{N}^+$ with $k_0< k_1<\cdots<k_m$. For all $0\le p\le h$ and $q\in\mathbb{N}$,  consider
\begin{align*}
{^{k_m}H_{pq}}=\bigcup_{p_{m-1}=0}^{h}
\bigcup_{p_{m-2}=0}^{h}\cdots
\bigcup_{p_{0}=0}^{h}
\bigcup_{\substack{
         q_y\in\mathbb{N}, \\
         1\le y\le m-1
         }} \, 
^{k_m}H_{pq} \cap{} ^{k_{m-1}}H_{p_{m-1}q_{m-1}} \cap{} \cdots \cap {}^{k_{0}}H_{p_{0}q_{0}}.
\end{align*}

Note that for all $1 \le j \le m$, if the intersection
\[
^{k_m}H_{pq} \cap\, ^{k_{m-1}}H_{p_{m-1}q_{m-1}} \cap \cdots \cap\, ^{k_j}H_{p_jq_j}
\]
is nonempty, then it admits the following decomposition:
\begin{align*}
&^{k_m}H_{pq} \cap\, ^{k_{m-1}}H_{p_{m-1}q_{m-1}} \cap \cdots \cap\, ^{k_j}H_{p_jq_j} \\
=~ &\bigcup_{p_{j-1} = 0}^{h} \bigcup_{q_{j-1} \in \mathbb{N}} 
\left( ^{k_m}H_{pq} \cap\, ^{k_{m-1}}H_{p_{m-1}q_{m-1}} \cap \cdots \cap\, ^{k_j}H_{p_jq_j} \right) \cap\, ^{k_{j-1}}H_{p_{j-1}q_{j-1}}.
\end{align*}

The number of nonempty sets in this decomposition is at most $D_H(k_j)$.  
If empty sets are allowed, we denote these sets by
\[
^{pq}Z_1^{\big\{ (p_{m-1},q_{m-1})\cdots (p_{j},q_{j}) \big\}},\ 
^{pq}Z_2^{\big\{ (p_{m-1},q_{m-1})\cdots (p_{j},q_{j}) \big\}},\cdots,\ 
^{pq}Z_{D_H(k_j)}^{\big\{ (p_{m-1},q_{m-1})\cdots (p_{j},q_{j}) \big\}}.
\]
This covers the case when $j \le m - 1$. If $j = m$, we simply write
\[
^{pq}Z_1,\ 
^{pq}Z_2,\ \cdots,\ 
^{pq}Z_{D_H(k_m)}.
\]

First, for any $p_{m-1}, \cdots, p_1 \in \{0, 1, \cdots, h\}$ and any $q_{m-1}, \cdots, q_1 \in \mathbb{N}$ such that the intersection $^{k_m}H_{pq} \cap\, ^{k_{m-1}}H_{p_{m-1}q_{m-1}} \cap \cdots \cap\, ^{k_1}H_{p_1q_1}$ is nonempty, we perform the following decomposition for each
\[
\bigg( \bigoplus_{I}G,~ ^{pq}Z_i^{\big\{ (p_{m-1},q_{m-1})\cdots (p_{1},q_{1}) \big\}} \bigg) \quad \text{with } 1 \le i \le D_H(k_1),
\]
then perform a saturated union operation to obtain some new families that cover the set $\big(\bigoplus_{I}G,~^{k_m}H_{pq} \cap\, ^{k_{m-1}}H_{p_{m-1}q_{m-1}} \cap \cdots \cap\, ^{k_1}H_{p_1q_1}\big)$.

According to Lemma~\ref{lemma: decomp for the infinite sum (pre for wr prod)}, for the given $k_0< k_1<\cdots<k_m$, there exist $J_i = J_i(k_0, k_{i-1})\in\mathbb{N}^+$ for all $1\le i\le m$, and $B_1(k_m)$-bounded families
\begin{align*}
\underbrace{
^{pq}\mathcal{U}_0^{\scriptstyle \substack{
    p_{2\sim m-1} \\
    q_{2\sim m-1}
    }}[{\scriptstyle(p_1,q_1,1)}]
}_{k_0\text{-disjoint}},\ 
\underbrace{
^{pq}\mathcal{U}_{r,j_r}^{\scriptstyle \substack{
    p_{2\sim m-1} \\
    q_{2\sim m-1}
    }}[{\scriptstyle(p_1,q_1,1)}]
}_{k_r\text{-disjoint}},\ 
\underbrace{
^{pq}\mathcal{U}^{\scriptstyle \substack{
    p_{2\sim m-1} \\
    q_{2\sim m-1}
    }}_{m,j_{m}}[{\scriptstyle(p_1,q_1,1)}]
}_{k_{m}\text{-disjoint}}
\end{align*}
which together cover $\bigg( \bigoplus_{I}G,~ ^{pq}Z_1^{\big\{ (p_{m-1},q_{m-1})\cdots (p_{1},q_{1}) \big\}} \bigg)$, where $1 \le r \le m-1$, and $1 \le j_i \le J_i$ for all $1 \le i \le m$.

Now consider a new sequence $k_0< k_1<\cdots<k_{m-1}<5B_1(k_m)$, and there exist $D_1(k_m)$-bounded families
\begin{align*}
\underbrace{
^{pq}\mathcal{U}_0^{\scriptstyle \substack{
    p_{2\sim m-1} \\
    q_{2\sim m-1}
    }}[{\scriptstyle(p_1,q_1,2)}]
}_{k_0\text{-disjoint}},\ 
\underbrace{
^{pq}\mathcal{U}_{r,j_r}^{\scriptstyle \substack{
    p_{2\sim m-1} \\
    q_{2\sim m-1}
    }}[{\scriptstyle(p_1,q_1,2)}]
}_{k_r\text{-disjoint}},\ 
\underbrace{
^{pq}\mathcal{U}^{\scriptstyle \substack{
    p_{2\sim m-1} \\
    q_{2\sim m-1}
    }}_{m,j_{m}}[{\scriptstyle(p_1,q_1,2)}]
}_{5B_1(k_{m})\text{-disjoint}}
\end{align*}
which cover $\bigg( \bigoplus_{I}G, ~ ^{pq}Z_2^{\big\{ (p_{m-1},q_{m-1})\cdots (p_{1},q_{1}) \big\}} \bigg)$, where $1 \le r \le m-1$, and $1 \le j_i \le J_i$ for all $1 \le i \le m$.
For each $1\le j_m\le J_m$, let 
\[
    ^{pq}\mathcal{V}^{\scriptstyle \substack{
    p_{2\sim m-1} \\
    q_{2\sim m-1}
    }}_{m,j_{m}}[{\scriptstyle(p_1,q_1,2)}]=
    ^{pq}\mathcal{U}^{\scriptstyle \substack{
    p_{2\sim m-1} \\
    q_{2\sim m-1}
    }}_{m,j_{m}}[{\scriptstyle(p_1,q_1,2)}]\cup_{k_m}{}
    ^{pq}\mathcal{U}^{\scriptstyle \substack{
    p_{2\sim m-1} \\
    q_{2\sim m-1}
    }}_{m,j_{m}}[{\scriptstyle(p_1,q_1,1)}],
\]
where $\cup_{k_m}$ denotes the $k_m$-saturated union.

Then, by Lemma~\ref{lemma: r-saturated}, the family $^{pq}\mathcal{V}^{\scriptstyle \substack{
    p_{2\sim m-1} \\
    q_{2\sim m-1}
    }}_{m,j_{m}}[{\scriptstyle(p_1,q_1,2)}]$ is $k_m$-disjoint, and we may assume it is $B_2(k_m)$-bounded. 

Similarly, consider a new sequence $k_0< k_1<\cdots<k_{m-1}<5B_2(k_m)$, and there exist $D_2(k_m)$-bounded families
\begin{align*}
\underbrace{
^{pq}\mathcal{U}_0^{\scriptstyle \substack{
    p_{2\sim m-1} \\
    q_{2\sim m-1}
    }}[{\scriptstyle(p_1,q_1,3)}]
}_{k_0\text{-disjoint}},\ 
\underbrace{
^{pq}\mathcal{U}_{r,j_r}^{\scriptstyle \substack{
    p_{2\sim m-1} \\
    q_{2\sim m-1}
    }}[{\scriptstyle(p_1,q_1,3)}]
}_{k_r\text{-disjoint}},\ 
\underbrace{
^{pq}\mathcal{U}^{\scriptstyle \substack{
    p_{2\sim m-1} \\
    q_{2\sim m-1}
    }}_{m,j_{m}}[{\scriptstyle(p_1,q_1,3)}]
}_{5B_2(k_{m})\text{-disjoint}}
\end{align*}
which cover $\bigg( \bigoplus_{I}G, ~ ^{pq}Z_3^{\big\{ (p_{m-1},q_{m-1})\cdots (p_{1},q_{1}) \big\}} \bigg)$, where $1 \le r \le m-1$, and $1 \le j_i \le J_i$ for all $1 \le i \le m$. 
For each $1\le j_m\le J_m$, let 
\[
    ^{pq}\mathcal{V}^{\scriptstyle \substack{
    p_{2\sim m-1} \\
    q_{2\sim m-1}
    }}_{m,j_{m}}[{\scriptstyle(p_1,q_1,3)}]=
    ^{pq}\mathcal{U}^{\scriptstyle \substack{
    p_{2\sim m-1} \\
    q_{2\sim m-1}
    }}_{m,j_{m}}[{\scriptstyle(p_1,q_1,3)}]\cup_{k_m}{}
    ^{pq}\mathcal{V}^{\scriptstyle \substack{
    p_{2\sim m-1} \\
    q_{2\sim m-1}
    }}_{m,j_{m}}[{\scriptstyle(p_1,q_1,2)}].
\]
Then, by Lemma~\ref{lemma: r-saturated}, the family $^{pq}\mathcal{V}^{\scriptstyle \substack{
    p_{2\sim m-1} \\
    q_{2\sim m-1}
    }}_{m,j_{m}}[{\scriptstyle(p_1,q_1,2)}]$ is $k_m$-disjoint, and we may assume it is $B_3(k_m)$-bounded. 

By analogy, we obtain $D_{D_H(k_1)}(k_m)+B_{D_H(k_1)}(k_m)$-bounded families 
\begin{align*}
\underbrace{
^{pq}\mathcal{U}_0^{\scriptstyle \substack{
    p_{2\sim m-1} \\
    q_{2\sim m-1}
    }}[{\scriptstyle(p_1,q_1,D_H(k_1))}]
}_{k_0\text{-disjoint}},\ 
\underbrace{
^{pq}\mathcal{U}_{r,j_r}^{\scriptstyle \substack{
    p_{2\sim m-1} \\
    q_{2\sim m-1}
    }}[{\scriptstyle(p_1,q_1,D_H(k_1))}]
}_{k_r\text{-disjoint}},\ 
\underbrace{
^{pq}\mathcal{V}^{\scriptstyle \substack{
    p_{2\sim m-1} \\
    q_{2\sim m-1}
    }}_{m,j_{m}}[{\scriptstyle(p_1,q_1,D_H(k_1))}]
}_{k_{m}\text{-disjoint}}
\end{align*}
which cover $\bigg( \bigoplus_{I}G, ~ ^{pq}Z_{D_H(k_1)}^{\big\{ (p_{m-1},q_{m-1})\cdots (p_{1},q_{1}) \big\}} \bigg)$, where $1 \le r \le m-1$, and $1 \le j_i \le J_i$ for all $1 \le i \le m$. Note that for each $1\le j_m\le J_m$, the family $^{pq}\mathcal{V}^{\scriptstyle \substack{
    p_{2\sim m-1} \\
    q_{2\sim m-1}
    }}_{m,j_{m}}[{\scriptstyle(p_1,q_1,D_H(k_1))}]$ is obtained by  
\[
    ^{pq}\mathcal{U}^{\scriptstyle \substack{
    p_{2\sim m-1} \\
    q_{2\sim m-1}
    }}_{m,j_{m}}[{\scriptstyle(p_1,q_1,D_H(k_1))}]\cup_{k_m}{}
    ^{pq}\mathcal{V}^{\scriptstyle \substack{
    p_{2\sim m-1} \\
    q_{2\sim m-1}
    }}_{m,j_{m}}[{\scriptstyle(p_1,q_1,D_H(k_1)-1)}].
\]

Furthermore, the families 
\begin{equation}\label{equ:covers}
\begin{aligned}
^{pq}\mathcal{U}_0^{\scriptstyle \substack{
    p_{2\sim m-1} \\
    q_{2\sim m-1}
    }} [{\scriptstyle(p_1, q_1, i_1)}],\ 
^{pq}\mathcal{U}_{r,j_r}^{\scriptstyle \substack{
    p_{2\sim m-1} \\
    q_{2\sim m-1}
    }}[{\scriptstyle(p_1, q_1, i_1)}],\ 
^{pq}\mathcal{V}^{\scriptstyle \substack{
    p_{2\sim m-1} \\
    q_{2\sim m-1}
    }}_{m,j_{m}}[{\scriptstyle(p_1,q_1,D_H(k_1))}]
\end{aligned}
\end{equation}
with ${1\le r\le m-1}$, $1\le j_r\le J_r$, $1\le j_m\le J_m$, and $1\le i_1\le D_H(k_1)$, together cover 
\[
\bigg( \bigoplus_{I}G, ~ ^{k_m}H_{pq} \cap\, ^{k_{m-1}}H_{p_{m-1}q_{m-1}} \cap \cdots \cap\, ^{k_1}H_{p_1q_1} \bigg).
\]

Note that the intersection $^{k_m}H_{pq} \cap\, ^{k_{m-1}}H_{p_{m-1}q_{m-1}} \cap \cdots \cap\, ^{k_1}H_{p_1q_1}$ coincides with the set $^{pq}Z_r^{\big\{ (p_{m-1}, q_{m-1}), \cdots, (p_2, q_2) \big\}}$ for some $1 \le r \le D_H(k_2)$. Moreover, for each such set, there exists a decomposition analogous to that in \cref{equ:covers}.

Therefore, by ranging over $0\le p_1 \le h$ and $q_1 \in \mathbb{N}$, we obtain the following uniformly bounded families: 
\begin{equation}
\begin{aligned}
\underbrace{
^{pq}\mathcal{U}_0^{\scriptstyle \substack{
    p_{3\sim m-1} \\
    q_{3\sim m-1}
    }} \Big[{\scriptstyle \substack{
    \cup(p_1, q_1, i_1) \\
    (p_2, q_2, i_2)
    }}\Big]
}_{k_0\text{-disjoint}},\quad 
\underbrace{
^{pq}\mathcal{U}_{r, j_r}^{\scriptstyle \substack{
    p_{3\sim m-1} \\
    q_{3\sim m-1}
    }}\Big[{\scriptstyle \substack{
    \cup(p_1, q_1, i_1) \\
    (p_2, q_2, i_2)
    }}\Big]
}_{k_r\text{-disjoint}},\quad 
\underbrace{
^{pq}\mathcal{V}_{m, j_m}^{\scriptstyle \substack{
    p_{3\sim m-1} \\
    q_{3\sim m-1}
    }} \Big[{\scriptstyle \substack{
    \cup(p_1, q_1, D_H(k_1)) \\
    (p_2, q_2, D_H(k_2))
    }}\Big] 
}_{k_m\text{-disjoint}}
\end{aligned}
\end{equation}
with $1 \le r \le m - 1$, $1 \le j_r \le J_r$, $1 \le j_m \le J_m$, $1 \le i_1 \le D_H(k_1)$, and $1 \le i_2 \le D_H(k_2)$.
These families together cover
\[
\bigg( \bigoplus_{I} G,\; ^{k_m}H_{pq} \cap\, ^{k_{m-1}}H_{p_{m-1}q_{m-1}} \cap \cdots \cap\, ^{k_2}H_{p_2q_2} \bigg).
\]

By analogy, by successively ranging over each pair $(p_i, q_i)$ with $p_i \in \{0,1,\cdots,h\}$ and $q_i \in \mathbb{N}$ for $2 \le i \le m - 1$, we obtain the following uniformly bounded families:
\begin{equation}
\scalebox{0.9}{$
\begin{aligned}
\underbrace{^{pq}\mathcal{U}_0\Bigg[{\scriptstyle \substack{
    \cup(p_1, q_1, i_1) \\
    \cup(p_2, q_2, i_2)\\
    \cdots \\
    \cup(p_{m-1}, q_{m-1}, i_{m-1})\\ i_m 
    }}\Bigg]}_{k_0\text{-disjoint}},\quad 
\underbrace{^{pq}\mathcal{U}_{r, j_r}\Bigg[{\scriptstyle \substack{
    \cup(p_1, q_1, i_1) \\
    \cup(p_2, q_2, i_2) \\
    \cdots \\
    \cup(p_{m-1}, q_{m-1}, i_{m-1})\\ i_m 
    }}\Bigg]}_{k_r\text{-disjoint}},\quad 
\underbrace{
^{pq}\mathcal{V}_{m, j_m}\Bigg[{\scriptstyle \substack{
    \cup(p_1, q_1, D_H(k_1)) \\
    \cup(p_2, q_2, D_H(k_2)) \\
    \cdots \\
    \cup(p_{m-1}, q_{m-1}, D_H(k_{m-1})) \\ 
    D_H(k_m)
    }}\Bigg] 
}_{k_m\text{-disjoint}}
\end{aligned}
$}
\end{equation}
with $1 \le r \le m - 1$, $1 \le j_r \le J_r$, $1 \le j_m \le J_m$, and $1 \le i_s \le D_H(k_s)$ for all $1\le s\le m$. 
These families together cover $\bigg( \bigoplus_{I} G,\; ^{k_m}H_{pq} \bigg)$. 

Finally, by ranging over $0\le p\le h$ and $q \in \mathbb{N}$, we obtain the following uniformly bounded family that covers $\bigg( \bigoplus_{I} G,\; H \bigg)$. 
\begin{enumerate}
    \item For each $0\le t \le h$, let
        \begin{align*}
        \mathcal{U}_{t} := \bigcup_{q \in \mathbb{N}}\,
        \bigcup_{\substack{
         0\le p_x\le h, \\
         1\le x\le m-1
         }} \,
        \bigcup_{\substack{
         q_y\in\mathbb{N}, \\
         1\le y\le m-1
         }} \,
        \bigcup_{i_1,\cdots,i_{m}}
        {^{tq}}\mathcal{U}_0\Bigg[{\scriptstyle \substack{
            \cup(p_1, q_1, i_1) \\
            \cup(p_2, q_2, i_2) \\
            \cdots \\
            \cup(p_{m-1}, q_{m-1}, i_{m-1})\\ i_m 
        }}\Bigg],
        \end{align*}
        where $1 \le i_s \le D_H(k_s)$ for all $1 \le s \le m$. 
        This yields $(h+1)$ many $k_0$-disjoint families. 
        Without loss of generality, we may assume that in the process of ranging over $p_1,\cdots,p_{m-1}$ and $q, q_1,\cdots,q_{m-1}$, we always consider those values for which the intersection $^{k_m}H_{tq} \cap\, ^{k_{m-1}}H_{p_{m-1}q_{m-1}} \cap \cdots \cap\, ^{k_1}H_{p_1q_1}$ is nonempty. The same assumption applies in the following context.

    \item For each $0\le v \le h$ and each $1 \le j_1 \le J_1$, let
        \begin{align*}
        \mathcal{U}_{1, j_1, v} := \bigcup_{q \in \mathbb{N}}\,
        \bigcup_{\substack{
         0\le p_x\le h, \\
         1\le x\le m-1
         }} \,
        \bigcup_{\substack{
         q_y\in\mathbb{N}, \\
         1\le y\le m-1
         }} \,
        \bigcup_{i_1,\cdots,i_{m}}
        {^{vq}}\mathcal{U}_{1, j_1}\Bigg[{\scriptstyle \substack{
            \cup(p_1, q_1, i_1) \\
            \cup(p_2, q_2, i_2) \\
            \cdots \\
            \cup(p_{m-1}, q_{m-1}, i_{m-1})\\ i_m 
        }}\Bigg],
        \end{align*}
        where $1 \le i_s \le D_H(k_s)$ for all $1 \le s \le m$.
        This yields $(h+1) \times J_1(k_0)$ many $k_1$-disjoint families.

    \item Similarly, for each $2 \le r \le m - 1$, we construct the following family:

        For each $0\le v \le h$, $1 \le j_r \le J_r$, and $1 \le i_w \le D_H(k_w)$ for $1 \le w \le r-1$, define
        \begin{align*}
        \mathcal{U}_{r, j_r, v, i_1, \cdots, i_{r-1}} := \bigcup_{q \in \mathbb{N}} \,
        \bigcup_{\substack{
         0\le p_x\le h, \\
         1\le x\le m-1
         }} \,
         \bigcup_{\substack{
         q_y\in\mathbb{N}, \\
         1\le y\le m-1
         }} \,
        \bigcup_{i_r,\cdots,i_{m}}
        {^{vq}}\mathcal{U}_{r, j_r}\Bigg[{\scriptstyle \substack{
            \cup(p_1, q_1, i_1) \\
            \cup(p_2, q_2, i_2) \\
            \cdots \\
            \cup(p_{m-1}, q_{m-1}, i_{m-1})\\ i_m 
        }}\Bigg],
        \end{align*}
        where $1 \le i_s \le D_H(k_s)$ for all $r \le s \le m$.
        This yields
        \[
        (h+1) \times J_r(k_0, k_{r-1}) \times \prod_{1 \le w \le r-1} |B_H(1_H, D(k_w))|
        \]
        many $k_r$-disjoint families.

    \item For each $0\le v \le h$ and $1 \le j_m \le J_m$, define
        \begin{align*}
        \mathcal{U}_{m, j_m, v} := 
        \bigcup_{q\in \mathbb{N}}\,
        \bigcup_{\substack{
         0\le p_x\le h, \\
         1\le x\le m-1
         }} \,
        \bigcup_{\substack{
         q_y\in\mathbb{N}, \\
         1\le y\le m-1
         }} \,
    ^{vq}\mathcal{V}_{m, j_m}\Bigg[{\scriptstyle \substack{
    \cup(p_1, q_1, D_H(k_1)) \\
    \cup(p_2, q_2, D_H(k_2)) \\
    \cdots \\
    \cup(p_{m-1}, q_{m-1}, D_H(k_{m-1})) \\ 
    D_H(k_m)
    }}\Bigg] 
        \end{align*}
        This yields $(h+1) \times J_m(k_0, k_{m-1})$ many $k_m$-disjoint families.
\end{enumerate}
Therefore, by \cref{lemma:more general decomp for trasdim}, we conclude that $\trasdim(W) \le \omega {m} + h$, and hence $W$ has asymptotic property C. 
\end{proof}

As a generalization of the result that the group $\mathbb{Z} \wr \mathbb{Z}$ has asymptotic property C, we show that the generalized conditions stated in the main theorem are valid, and provide examples to illustrate this. 

\begin{example}[FC-groups with finite asymptotic dimension]
According to \cref{lem:fgFChasFAD}, every finitely generated FC-group has finite asymptotic dimension. This class includes, for instance, all finite groups and all finitely generated virtually abelian groups. If the assumption of ``finitely generated'' is dropped, then the direct sum of countably many groups of asymptotic dimension zero (equipped with a proper left-invariant metric) still has finite asymptotic dimension.

Note that the class of FC-groups is closed under operations such as taking subgroups, direct products, and countable direct sums. 
\end{example}

\begin{example}[Wreath-like product with a symmetric transversal and a unital, bornologous transversal map]

In our main theorem, we require that the action in $W \in \WR(G, H \curvearrowright I)$ has finitely many orbits, and that there exists a symmetric transversal $T$ associated with ${G}^{(I)} \lhd W$, such that the transversal map $\eta: H \longrightarrow T$ is unital and bornologous. This class of groups $W$ naturally generalizes the classical $G\wr H$. 

As an example, consider the following non-split exact sequence:
\begin{align*}
1\longrightarrow 2\mathbb{Z}\longrightarrow \mathbb{Z}\longrightarrow \mathbb{Z}_2\longrightarrow 1.
\end{align*}
It naturally extends to a wreath-like product:
\begin{align*}
1\longrightarrow \bigoplus_{\mathbb{Z}\cup\{*\}}\mathbb{Z}\longrightarrow &(\mathbb{Z}\wr\mathbb{Z})\oplus \mathbb{Z} \longrightarrow \mathbb{Z}\oplus \mathbb{Z}_2\longrightarrow 1,\\
\big((a_i)_{i\in\mathbb{Z}},p\big)\mapsto &\big[ \big((a_i)_{i\in\mathbb{Z}},0\big),2p \big],\\ 
&\big[ \big((a_i)_{i\in\mathbb{Z}},b\big),h \big]\mapsto (b, h\text{ mod } 2).
\end{align*}

The action of $\mathbb{Z} \oplus \mathbb{Z}_2$ on the set $\mathbb{Z} \cup \{*\}$ is defined as follows: for all $(n,h) \in \mathbb{Z} \oplus \mathbb{Z}_2$ and $m \in \mathbb{Z}$, set $(n,h). m = n + m$ and $(n,h). * = *$.

Clearly, $(\mathbb{Z}\wr\mathbb{Z})\oplus \mathbb{Z}$ is not a classical wreath product of $\mathbb{Z}$ with $\mathbb{Z} \oplus \mathbb{Z}_2$.
Let $\vec{0}$ denote the unit element of $\bigoplus_{\mathbb{Z}} \mathbb{Z}$. Define a transversal $T:=\big\{ \big[ (\vec{0},b),h \big] \mid  b\in \mathbb{Z}, h\in \mathbb{Z}_2 \big\}$ and a transversal map $\eta: \mathbb{Z} \oplus \mathbb{Z}_2\longrightarrow T$ by
$(b,h)\mapsto \big[ (\vec{0},b),h \big]$. 
Then it is easy to check that $\eta$ is bornologous. 
\end{example}

The following section shows the examples of groups with the $(m,0)$-DTUT property. 

\section{Groups with the DTUT Property}\label{sec:Groups with DTUT}
Since uniform disjointness is a coarse invariant, similarly to the permanence properties of finite asymptotic dimension, we can directly obtain the following propositions.

\begin{proposition}\label{prop:(m,n)-DTUT is a coarse inv}
For any $m,n\in\mathbb{N}$, the $(m,n)$-DTUT property is a coarse invariant. 
\end{proposition}

\begin{proposition}\label{prop:(m,n)-DTUT is closed under product}
Let $\Gamma_1,\Gamma_2$ be two groups with $(m,n_1)$-DTUT and $(m,n_2)$-DTUT, respectively, where $m,n_1,n_2\in\mathbb{N}^+$, then $\Gamma_1\times \Gamma_2$ has the $(m,n_1\times n_2+n_1+n_2)$-DTUT property. 
\end{proposition}

We now give some examples of groups satisfying the $(m,n)$-DTUT property. We begin by verifying that $\mathbb{Z}^n$ (for $n \in \mathbb{N}^+$) satisfies the $(1,0)$-DTUT property.

\begin{example}\label{eg:Z^n has (1,0)-DTUT}
For all $n \in \mathbb{N}^+$, it is easy to verify that $\mathbb{Z}^n$ satisfies the $(1,0)$-DTUT property.

Let $h \colon \mathbb{N}^+ \to \mathbb{N}^+$ be an arbitrary map.  
For every pair $k_0 < k_1$ with $k_0, k_1 \in \mathbb{N}^+$, and for every $1\le l \le h(k_1)$, define
\begin{align*}
&x_l:=\Big(2k_1(l-1),2k_1(l-1),\cdots,2k_1(l-1)\Big)\in\mathbb{Z}^n,\\
&S:=(2{h(k_1)}+1) k_1\in\mathbb{Z}.    
\end{align*}
Then, we define 
\begin{align*}
\mathcal{C}^l:=\bigg\{
\prod_{j=1}^n 
\Big[{S} i_j+k_1,~{S}(i_j+1)-k_1 \Big)+x_l ~\bigg|~ i_j\in\mathbb{Z}\bigg\}.
\end{align*}   

Note that every element in $\mathcal{C}_l$ is an $n$-hypercube. 
If $n=1$, let 

\begin{align*}
\mathcal{D}_{1,1,0}^l:=\bigg\{
& \Big[{S} i_1, ~{S} i_1+k_1 \Big)+x_l ~\bigg|~  i_1\in\mathbb{Z}\bigg\},\\
\mathcal{D}_{1,1,1}^l:=\bigg\{
& \Big[{S} (i_1{+}1)-k_1, ~{S} (i_1{+}1) \Big)+x_l ~\bigg|~  i_1\in\mathbb{Z}\bigg\}.
\end{align*}

If $n\ge 2$, consider a bijection $\psi:\ \{1,2,\cdots,2^{n-1}\}\longrightarrow \{0,1\}^{n-1}$, and for all $1\le u \le 2^{n-1}$, $1\le p\le n-1$, and $r\in\mathbb{Z}$, define
\begin{align*}
T_u(r,p):= 2r+\psi(u)_{p}\in\mathbb{Z}.
\end{align*}

Then, for all $u\in\{1,2,\cdots,2^{n-1}\}$, we define 

\begin{align*}
\mathcal{D}_{u,1,0}^l:=\bigg\{\bigg[
& \Big[{S} i_1, ~{S} i_1+k_1 \Big)\times \prod_{2\le q \le n}
{S}\Big[T_u(i_q, q-1), ~T_u(i_q, q-1)+1 \Big)
\bigg]+x_l ~\bigg|~ i_1,i_q\in\mathbb{Z} \bigg\},\\
\mathcal{D}_{u,1,1}^l:=\bigg\{\bigg[
& \Big[{S} (i_1{+}1)-k_1, ~{S} (i_1{+}1) \Big)\times\\
&\prod_{2\le q \le n}
{S}\Big[T_u(i_q, q-1), ~T_u(i_q, q-1)+1 \Big)
\bigg]+x_l ~\bigg|~ i_1,i_q\in\mathbb{Z} \bigg\}.\\
\end{align*}   
Furthermore, for all $v\in\{2,\cdots,n\}$, define 
\begin{align*}
\mathcal{D}_{u,v,0}^l:=\bigg\{\bigg[
&\prod_{\substack{j<v\\1\le j\le n}}
{S}\Big[T_u(i_j,j), T_u(i_j,j)+1 \Big)
\times \Big[{S} i_v,\ {S} i_v+k_1 \Big)\\
& \times \prod_{\substack{q>v\\1\le q \le n}}
{S}\Big[T_u(i_q, q-1), T_u(i_q, q-1)+1 \Big)
\bigg]+x_l   ~\bigg|~ i_j,i_v,i_q\in\mathbb{Z} \bigg\},\\
\mathcal{D}_{u,v,1}^l:=\bigg\{\bigg[
&\prod_{\substack{j<v\\1\le j\le n}}
{S}\Big[T_u(i_j,j),\ T_u(i_j,j)+1 \Big)
\times \Big[{S} (i_v+1)-k_1,~ {S} (i_v+1) \Big)\\
&\times \prod_{\substack{q>v\\1\le q \le n}}
{S}\Big[T_u(i_q, q-1), T_u(i_q, q-1)+1 \Big)
\bigg]+x_l ~\bigg|~ i_j,i_v,i_q\in\mathbb{Z} \bigg\}.
\end{align*}

Figure~\ref{fig:Z^2 has (1,0)-DTUT} shows an example of the $(1,0)$-DTUT decomposition of $\mathbb{Z}^2$. For a fixed $l$, the interior of each small square, represented by the black shaded region, corresponds to a subset of the family $\mathcal{C}^l$.

Each shell surrounding the shaded region is composed of two pairs of colored strip regions. By selecting a coordinate axis direction, one determines the pair of strips that form the shell (the $x$-axis corresponds to the green and purple strips, while the $y$-axis corresponds to the yellow and red pair). This choice of axis corresponds to the index $v$ (ranging from $1$ to $n$ in the case of $\mathbb{Z}^n$, and the shell of an $n$-dimensional cube consists of $n$ such pairs).

The parameter $z \in \{0,1\}$ specifies which color in each pair is selected. The value of $u$ (ranging from $1$ to $2^{n-1}$ in the case of $\mathbb{Z}^n$) encodes the parity configuration along the remaining $n-1$ coordinate directions, which determines the positions of the colored strips in those axes.

With the above understanding, it is easy to verify that 
\begin{itemize}
  \item $\forall~1\le l\le {h(k_1)}$, the family $\bigg( \bigcup_{u=1}^{2^{n-1}}\bigcup_{v=1}^n\bigcup_{z\in\{0,1\}}\mathcal{D}_{l,u,v,z}\bigg)\cup \mathcal{C}_l$ covers $\mathbb{Z}^n$. 
  \item $\forall~1\le l\le {h(k_1)}$, $\mathcal{C}_l$ is $2k_1$-disjoint, and $2^{n-1}\times (2(n+1)h(k_1)+1) k_1$-bounded.
  \item $\forall~1\le u \le 2^{n-1}$, $1\le v\le n$, and $z\in\{0,1\}$, the family $\bigcup_{l=1}^{{n+1}^{|T(m,k)|}}\mathcal{D}_{l,u,v,z}$ is $k_1$-disjoint, and $2^{n-1}\times (2(n+1)h(k_1)+1)k_1$-bounded.
\end{itemize}   
\end{example}

\begin{figure}[h]
\centering
\begin{tikzpicture}[scale=0.8]

\def\cellsize{2}
\def\squaresize{1.2}
\def\offset{0.4}
\def\ext{1}

\foreach \i in {0,1,2,3} {
  \draw[thick] (\i*\cellsize,-\ext) -- (\i*\cellsize,3*\cellsize+\ext);
  \draw[thick] (-\ext,\i*\cellsize) -- (3*\cellsize+\ext,\i*\cellsize);
}

\foreach \x/\y in {0/0, 0/2, 2/0, 2/2} {
  \begin{scope}[shift={(\x*\cellsize,\y*\cellsize)}]
    \fill[red, opacity=0.3] (0,\offset+\squaresize) rectangle (\cellsize,\cellsize);
    \fill[yellow, opacity=0.3] (0,0) rectangle (\cellsize,\offset);
    
    \fill[green!60!black, opacity=0.3] (0,0) rectangle (\offset,\cellsize);
    \fill[blue!60!black, opacity=0.3] (\offset+\squaresize,0) rectangle (\cellsize,\cellsize);
    
    \fill[pattern=north east lines, pattern color=black] 
      (\offset,\offset) rectangle (\offset+\squaresize,\offset+\squaresize);
    \draw[thick] (\offset,\offset) rectangle (\offset+\squaresize,\offset+\squaresize);
  \end{scope}
}

\foreach \x/\y in {1/0, 0/1, 1/1, 2/1, 1/2} {
  \begin{scope}[shift={(\x*\cellsize,\y*\cellsize)}]    
    \fill[pattern=north east lines, pattern color=black] 
      (\offset,\offset) rectangle (\offset+\squaresize,\offset+\squaresize);
    \draw[thick] (\offset,\offset) rectangle (\offset+\squaresize,\offset+\squaresize);
  \end{scope}
}

\foreach \x/\y in {0/1, 0/3, 2/1, 2/3} {
  \begin{scope}[shift={(\x*\cellsize,\y*\cellsize)}]
    \fill[yellow, opacity=0.3] (0,0) rectangle (\cellsize,\offset);
  \end{scope}
}
\foreach \x/\y in {0/-1, 0/1, 2/-1, 2/1} {
  \begin{scope}[shift={(\x*\cellsize,\y*\cellsize)}]
    \fill[red, opacity=0.3] (0,\offset+\squaresize) rectangle (\cellsize,\cellsize);
  \end{scope}
}
\foreach \x/\y in {-1/0, 1/0, 1/2, -1/2} {
  \begin{scope}[shift={(\x*\cellsize,\y*\cellsize)}]
    \fill[blue!60!black, opacity=0.3] (\offset+\squaresize,0) rectangle (\cellsize,\cellsize);
  \end{scope}
}
\foreach \x/\y in {1/0, 1/2, 3/0, 3/2} {
  \begin{scope}[shift={(\x*\cellsize,\y*\cellsize)}]
    \fill[green!60!black, opacity=0.3] (0,0) rectangle (\offset,\cellsize);
  \end{scope}
}

\draw[->, thick] 
  (0, \cellsize) 
  -- 
  (0, \cellsize*3.5) node[pos=1, left] {$y$};
\draw[->, thick] 
  (0, \cellsize) 
  -- 
  (\cellsize*2.2, \cellsize*3.2);
\draw[thick] 
  (0, \cellsize) 
  -- 
  (\cellsize*2.4, \cellsize*3.4);
\draw[->, thick] 
  (0, \cellsize) 
  -- 
  (\cellsize*3.5, \cellsize) node[pos=1, below] {$x$};
  
\end{tikzpicture}
\caption{$\mathbb{Z}^2$ has the $(1,0)$-DTUT property.}
\label{fig:Z^2 has (1,0)-DTUT}
\end{figure}
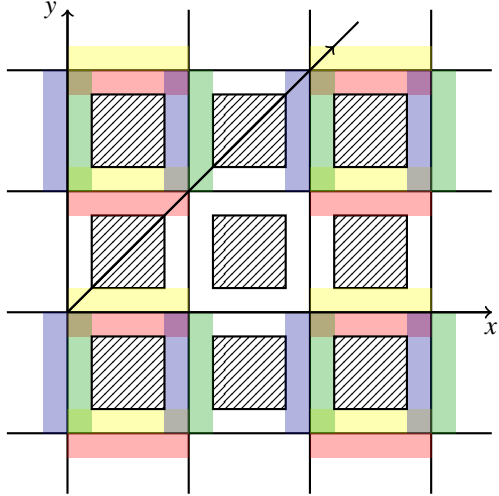

\subsection{Virtually Nilpotent Groups}
Recall that a group $G$ is called \emph{nilpotent} if it has a \emph{central series} of finite length. That is, a series of normal subgroups
\[
\{1\} = G_0 \triangleleft G_1 \triangleleft \cdots \triangleleft G_n = G,
\]
where $G_{i+1}/G_i \le Z(G/G_i)$, or equivalently, $[G, G_{i+1}] \le G_i$. Moreover, the smallest $n$ such that $G$ has a central series of length $n$ is called the \emph{nilpotency class} of $G$, and $G$ is said to be \emph{nilpotent of class $n$}. 

Since $\mathbb{Z}^n$ (for $n \in \mathbb{N}^+$) is a specific example of a nilpotent group, it is natural to ask whether every nilpotent group satisfies the $(m,0)$-DTUT property for some $m \in \mathbb{N}^+$. As a next step, we consider the case of the integer Heisenberg group $H_3(\mathbb{Z})$.

Recall that
\[
H_3(\mathbb{Z}) = \Bigg\{ 
\scalebox{0.8}{$
\begin{pmatrix}
1 & x & z \\
0 & 1 & y \\
0 & 0 & 1
\end{pmatrix}
$}
~\big|~ x, y, z \in \mathbb{Z}
\Bigg\}.
\]

Consider $a=
\scalebox{0.8}{$
\begin{pmatrix}
1 & 1 & 0 \\
0 & 1 & 0 \\
0 & 0 & 1
\end{pmatrix}
$}$, 
$b=
\scalebox{0.8}{$
\begin{pmatrix}
1 & 0 & 0 \\
0 & 1 & 1 \\
0 & 0 & 1
\end{pmatrix}
$}$, and 
$c=\scalebox{0.8}{$
\begin{pmatrix}
1 & 0 & 1 \\
0 & 1 & 0 \\
0 & 0 & 1
\end{pmatrix}
$}$.
Since $c = aba^{-1}b^{-1}$, the elements $a$ and $b$ are sufficient to generate $H_3(\mathbb{Z})$ and in fact yield the presentation
\[
H_3(\mathbb{Z}) = \langle a, b \mid b^{-1}aba^{-1} = aba^{-1}b^{-1} = ba^{-1}b^{-1}a \rangle.
\]

For convenience, we include the element $c$ as an additional (third) generator, resulting in the more commonly used presentation:
\[
H_3(\mathbb{Z}) = \langle a, b, c \mid aba^{-1}b^{-1} = c, \ ac = ca, \ bc = cb \rangle.
\]

Note that there exists a one-to-one correspondence between $H_3(\mathbb{Z})$ and $\mathbb{Z}^3$. In particular, when each coset of $ \langle z \rangle $ is shrunk to a single point, the resulting Cayley graph becomes isomorphic to that of $ \mathbb{Z}^2 $. 

\begin{proposition}\label{prop:H_3(Z) has (2,0)DTUT}
The integer Heisenberg group $H_3(\mathbb{Z})$ has the $(2,0)$-DTUT property. 
\end{proposition}
\begin{proof}
Note that $H_3(\mathbb{Z})$ admits the following central series:
\[
1 \triangleleft \langle c \rangle \triangleleft H_3(\mathbb{Z}),
\]
in particular, $\langle c \rangle \cong \mathbb{Z}$ and $H_3(\mathbb{Z}) / \langle c \rangle \cong \mathbb{Z}^2$. Moreover, both $\mathbb{Z}$ and $\mathbb{Z}^2$ satisfy the $(1,0)$-DTUT property.
Now let $h\colon \mathbb{N}^+ \to \mathbb{N}^+$ be an arbitrary map, and let $k_0 < k_1 < k_2$ be arbitrary positive integers.

Since $\mathbb{Z}^2$ satisfies the $(1,0)$-DTUT property, we consider $k_1 < k_2$, then for any $1\le l_2 \le h(k_2)$, there exists a $k_1$-disjoint family $\mathcal{C}_{\mathbb{Z}^2}^{l_2}$, and $k_2$-disjoint families
\[
\mathcal{D}_{\mathbb{Z}^2,0}^{l_2},\ \mathcal{D}_{\mathbb{Z}^2,1}^{l_2},\ \mathcal{D}_{\mathbb{Z}^2,2}^{l_2},\ \cdots,\ \mathcal{D}_{\mathbb{Z}^2,7}^{l_2},
\]
such that these families cover $\mathbb{Z}^2$, and for each $0 \leq i \leq 7$, the family $\bigcup_{l_2 = 1}^{h(k_2)} \mathcal{D}_{\mathbb{Z}^2,i}^{l_2}$ is $k_2$-disjoint.
We further assume that each of the above families is $M_{\mathbb{Z}^2}(k_2)$-bounded.

For each set $C^{l_2} \in \mathcal{C}_{\mathbb{Z}^2}^{l_2}$, fix a constant $y_{C^{l_2}} \in \mathbb{Z}$ such that for every $(x, y) \in C^{l_2}$, it holds that $y_{C^{l_2}} - M_{\mathbb{Z}^2}(k_1) \le y \le y_{C^{l_2}} + M_{\mathbb{Z}^2}(k_1)$. 
Similarly, for each set $D_i^{l_2} \in \mathcal{D}_{\mathbb{Z}^2,i}^{l_2}$, fix a constant $y_{D_i^{l_2}} \in \mathbb{Z}$ such that for all $(x, y) \in D_i^{l_2}$, we have $y_{D_i^{l_2}} - M_{\mathbb{Z}^2}(k_2) \le y \le y_{D_i^{l_2}} + M_{\mathbb{Z}^2}(k_2)$ for each $0 \le i \le 7$. 

Now for each fixed $l_2$, define a map $\psi^{l_2} \colon \mathbb{Z}^2 \to H_3(\mathbb{Z})$ as follows: if $(x, y) \in C^{l_2}$ for some $C^{l_2} \in \mathcal{C}_{\mathbb{Z}^2}^{l_2}$, then define 
\[
\psi^{l_2}(x, y) := a^x b^y c^{-x y_{C^{l_2}}}.
\]
Similarly, if $(x, y) \in D_i^{l_2}$ for some $D_i^{l_2} \in \mathcal{D}_{\mathbb{Z}^2,i}^{l_2}$, then define $\psi^{l_2}(x, y) := a^x b^y c^{-x  y_{D_i^{l_2}}}$.

It is straightforward to verify that $\psi^{l_2}(\mathcal{C}_{\mathbb{Z}^2}^{l_2})$ is $k_1$-disjoint, and $\big[ M_{\mathbb{Z}^2}(k_2) \big]^2 + 2 M_{\mathbb{Z}^2}(k_2)$-bounded in $H_3(\mathbb{Z})$.  
Indeed, for any $(x, y) \in C^{l_2}$ and $(x', y') \in \widetilde{C}^{l_2}$ with $C^{l_2} \ne \widetilde{C}^{l_2} \in \mathcal{C}_{\mathbb{Z}^2}^{l_2}$, we compute
\[
\big( \psi^{l_2}(x, y) \big)^{-1} \psi^{l_2}(x', y') 
= a^{x' - x} b^{y' - y} c^{x  y_{C^{l_2}} - x'  y_{\widetilde{C}^{l_2}} + y (x' - x)}.
\]
Since $\mathcal{C}_{\mathbb{Z}^2}^{l_2}$ is $k_1$-disjoint, it follows that $|x' - x| + |y' - y| \ge k_1$, which implies that the total number of occurrences of the word $a^{\pm1}$ or $b^{\pm1}$ in the shortest string is greater than $k_1$. Hence, $\psi^{l_2}(\mathcal{C}_{\mathbb{Z}^2}^{l_2})$ is $k_1$-disjoint.

On the other hand, if $(x, y), (x', y') \in C^{l_2}$ for some $ C^{l_2} \in \mathcal{C}_{\mathbb{Z}^2}^{l_2}$, then
\[
\big( \psi^{l_2}(x, y) \big)^{-1} \psi^{l_2}(x', y') 
= a^{x' - x} b^{y' - y} c^{(x - x')y_{C^{l_2}} + y (x' - x)}.
\]
We then have the estimate
\[
|x' - x| + |y' - y| + |y - y_{C^{l_2}}| \times |x' - x|
\le \big[ M_{\mathbb{Z}^2}(k_2) \big]^2 + 2 M_{\mathbb{Z}^2}(k_2).
\]
This shows that $\psi^{l_2}(\mathcal{C}_{\mathbb{Z}^2}^{l_2})$ is $\big[ M_{\mathbb{Z}^2}(k_2) \big]^2 + 2 M_{\mathbb{Z}^2}(k_2)$-bounded.

Similarly, for each $0 \le i \le 7$, the image $\psi^{l_2}(\mathcal{D}_{\mathbb{Z}^2,i}^{l_2})$ is $k_2$-disjoint and $\big[ M_{\mathbb{Z}^2}(k_2) \big]^2 + 2 M_{\mathbb{Z}^2}(k_2)$-bounded in $H_3(\mathbb{Z})$. Moreover, the family $\bigcup_{l_2 = 1}^{h(k_2)} \psi^{l_2} \big( \mathcal{D}_{\mathbb{Z}^2,i}^{l_2} \big)$ is $k_2$-disjoint. 

In particular, the map $\psi^{l_2}$ defines a ``local'' coarse embedding from $\mathbb{Z}^2$ into $H_3(\mathbb{Z})$, in the sense that it preserves both disjointness and uniform boundedness of families under the image.

Since $\mathbb{Z}$ has the $(1,0)$-DTUT property, consider $10 \big[ M_{\mathbb{Z}^2}(k_2) \big]^2 < 10 \big[ M_{\mathbb{Z}^2}(k_2) \big]^2 + 1$, then for any $1\le l_1\le h(k_1)$, there exists a $10 \big[ M_{\mathbb{Z}^2}(k_2) \big]^2$-disjoint family $\mathcal{C}_{\mathbb{Z}}^{l_1}$, as well as two $10 \big[ M_{\mathbb{Z}^2}(k_2) \big]^2$-disjoint families $\mathcal{D}_{\mathbb{Z},0}^{l_1},\ \mathcal{D}_{\mathbb{Z},1}^{l_1}$ such that these families cover $\mathbb{Z}$ and, for each $i = 0,1$, the family $\bigcup_{l_1 = 1}^{h(k_1)} \mathcal{D}_{\mathbb{Z},i}^{l_1}$ is $10 \big[ M_{\mathbb{Z}^2}(k_2) \big]^2$-disjoint in $\mathbb{Z}$. 

Now consider the map
\[
\phi \colon \mathbb{Z} \longrightarrow H_3(\mathbb{Z}), \quad z \longmapsto c^z.
\]
Note that $\phi$ is a coarse embedding, it is easy to calculate that $\phi\big( \mathcal{C}_{\mathbb{Z}}^{l_1} \big),\ \phi\big( \mathcal{D}_{\mathbb{Z},0}^{l_1} \big),\ \phi\big( \mathcal{D}_{\mathbb{Z},1}^{l_1} \big)$ are $3 M_{\mathbb{Z}^2}(k_2)$-disjoint and uniformly bounded in $H_3(\mathbb{Z})$. Let us denote this uniform bound by $N_{\mathbb{Z}}(k_2)$. 

For every $l_1 \in \{1, 2, \cdots, h(k_1)\}$ and $l_2 \in \{1, 2, \cdots, h(k_2)\}$, define
\begin{align*}
\mathcal{C}_0^{l_1,l_2} 
&:= \bigcup_{C^{l_2} \in \mathcal{C}_{\mathbb{Z}^2}^{l_2}} \psi^{l_2}(C^{l_2}) \cdot \phi\big( \mathcal{C}_{\mathbb{Z}}^{l_1} \big), \\
\mathcal{D}_{1,i}^{l_1,l_2} 
&:= \bigcup_{C^{l_2} \in \mathcal{C}_{\mathbb{Z}^2}^{l_2}} \psi^{l_2}(C^{l_2}) \cdot \phi\big( \mathcal{D}_{\mathbb{Z},i}^{l_1} \big), \quad (i \in \{0,1\}), \\
\mathcal{D}_{2,j}^{l_1,l_2} 
&:= \bigcup_{D_j^{l_2} \in \mathcal{D}_{\mathbb{Z}^2,j}^{l_2}} \psi^{l_2}(D_j^{l_2}) \cdot \phi\big( \mathcal{C}_{\mathbb{Z}}^{l_1} \big), \quad (j \in \{0,1,\cdots,7\}), \\
\mathcal{D}_{2,(u,v)}^{l_1,l_2} 
&:= \bigcup_{D_v^{l_2} \in \mathcal{D}_{\mathbb{Z}^2,v}^{l_2}} \psi^{l_2}(D_v^{l_2}) \cdot \phi\big( \mathcal{D}_{\mathbb{Z},u}^{l_1} \big), \quad (u \in \{0,1\},\ v \in \{0,1,\cdots,7\}).
\end{align*}

It is also straightforward to verify that $\mathcal{C}_0^{l_1,l_2}$ is $k_0$-disjoint, each $\mathcal{D}_{1,i}^{l_1,l_2}$ is $k_1$-disjoint for $i \in \{0,1\}$, and both $\mathcal{D}_{2,j}^{l_1,l_2}$ and $\mathcal{D}_{2,(u,v)}^{l_1,l_2}$ are $k_2$-disjoint for all $j \in \{0,1,\cdots,7\}$, $u \in \{0,1\}$, and $v \in \{0,1,\cdots,7\}$. Moreover, for fixed $l_2$, the family $\bigcup_{l_1 = 1}^{h(k_1)} \mathcal{D}_{1,i}^{l_1,l_2}$ is $k_1$-disjoint in $H_3(\mathbb{Z})$, and for fixed $l_1$, both $\bigcup_{l_2 = 1}^{h(k_2)} \mathcal{D}_{2,j}^{l_1,l_2}$ and $\bigcup_{l_2 = 1}^{h(k_2)} \mathcal{D}_{2,(u,v)}^{l_1,l_2}$ are $k_2$-disjoint in $H_3(\mathbb{Z})$. Furthermore, all of these families are uniformly bounded by
\[
\big[ M_{\mathbb{Z}^2}(k_2) \big]^2 + 2 M_{\mathbb{Z}^2}(k_2) + N_{\mathbb{Z}}(k_2).
\]

Indeed, the uniform boundedness is clear. We now verify that for fixed $l_2$, the family 
$\bigcup_{l_1 = 1}^{h(k_1)} \mathcal{D}_{1,i}^{l_1,l_2}$ is $k_1$-disjoint.
Note that $\bigcup_{l_1 = 1}^{h(k_1)} \mathcal{D}_{\mathbb{Z},i}^{l_1}$ is $10 \big[ M_{\mathbb{Z}^2}(k_2) \big]^2$-disjoint in $\mathbb{Z}$.

Let $C^{l_2}, \widetilde{C}^{l_2} \in \mathcal{C}_{\mathbb{Z}^2}^{l_2}$ and $F, \widetilde{F}\in \bigcup_{l_1 = 1}^{h(k_1)} \mathcal{D}_{\mathbb{Z},i}^{l_1}$.  
Choosing arbitrary points $(x,y) \in C^{l_2}$, $(x',y') \in \widetilde{C}^{l_2}$, $z \in F$, and $z' \in \widetilde{F}$, we consider the following cases: 
\begin{itemize}
    \item If $C^{l_2} \neq \widetilde{C}^{l_2}$, since $\psi^{l_2}(\mathcal{C}_{\mathbb{Z}^2}^{l_2})$ is $k_1$-disjoint, the shortest string $\big( \psi^{l_2}(x, y) \cdot \phi(z) \big)^{-1} \cdot \psi^{l_2}(x', y') \cdot \phi(z')$ must contain more than $k_1$ words of the form $a^{\pm 1}$ or $b^{\pm 1}$. Hence,
    \[
    d_{H_3(\mathbb{Z})}\big( \psi^{l_2}(C^{l_2}) \cdot \phi(F),~\psi^{l_2}(\widetilde{C}^{l_2}) \cdot \phi(\widetilde{F}) \big) > k_1.
    \]
    
    \item If $C^{l_2} = \widetilde{C}^{l_2}$, and $F \ne \widetilde{F}$, then since $d_{\mathbb{Z}}(F, \widetilde{F})> 10 \big[ M_{\mathbb{Z}^2}(k_2) \big]^2$, we have $|z' - z| > 10 \big[ M_{\mathbb{Z}^2}(k_2) \big]^2$. Note that 
    \begin{align*}
     &u := \psi^{l_2}(x, y) \cdot \phi(z) = a^x b^y c^{-x y_{C^{l_2}} + z}, \\
     &v := \psi^{l_2}(x', y') \cdot \phi(z') = a^{x'} b^{y'} c^{-x' y_{C^{l_2}} + z'}.   
    \end{align*}
    Then, $u^{-1} v = a^{x' - x} b^{y' - y} c^{(y - y_{C^{l_2}})(x' - x) + z' - z}$. Note that $|x' - x|, |y' - y| \le M_{\mathbb{Z}^2}(k_2)$, and $\big| z' - z + (y - y_{C^{l_2}})(x' - x) \big| > 9 \big[ M_{\mathbb{Z}^2}(k_2) \big]^2$. Therefore, we have
    \[
    l_{H_3(\mathbb{Z})}(c^{(y - y_{C^{l_2}})(x' - x) + z' - z}) > 3 M_{\mathbb{Z}^2}(k_2), \text{ and }
    l_{H_3(\mathbb{Z})}(u^{-1} v) > M_{\mathbb{Z}^2}(k_2) \ge k_2.
    \]
\end{itemize}
The disjointness of the remaining families can be verified in a similar manner.
\end{proof}

This result naturally extends to all finitely generated torsion-free nilpotent groups. 
\begin{theorem}\label{thm:f.g. torsion-free group has DTUT}
Let $G$ be a finitely generated torsion-free nilpotent group of class $m$, then $G$ has the $(m,0)$-DTUT property.   
\end{theorem}

\begin{proof}
We proceed by induction on $m$. When $m = 1$, the group $G \cong \mathbb{Z}^k$ for some $k \in \mathbb{N}^+$, and thus $G$ has the $(1,0)$-DTUT property.
Assume that every nilpotent group of class $m-1$ satisfies the $(m-1,0)$-DTUT property. Now let $G$ be a nilpotent group of class $m$. Then $G/Z(G)$ is nilpotent of class $m - 1$. The following table summarizes the central series, generators, and relations for $G$ and $G/Z(G)$:
{
\small
\[
\setlength{\arraycolsep}{4pt}
\renewcommand{\arraystretch}{1.0}
\begin{array}{|c|c|c|}
\hline
\makecell[c]{\rule{0pt}{3.8ex}\textbf{Central}\\\textbf{series}}
&
\displaystyle \{1\} = G_0 \triangleleft 
\underset{\substack{\rotatebox{90}{$=$} \\[-1pt] Z(G)}}{G_1}
\triangleleft \cdots \triangleleft G_m = G 
&
\displaystyle \{1\} = G_1/G_1 \triangleleft G_2/G_1 \triangleleft \cdots \triangleleft \underset{\substack{\rotatebox{90}{$=$} \\[-1pt] G/Z(G)}}{G_m/G_1} \\
\hline
\textbf{Generators} &
\begin{array}{c}
G_j = \langle s_1, \cdots, s_{i_j} \rangle\ (1 \le j \le m) \\
1 \le i_1 < \cdots < i_m = n
\end{array} &
\begin{array}{c}
G_t/G_1 = \langle \overline{s_{i_1{+}1}}, \cdots, \overline{s_{i_t}} \rangle\ (2 \le t \le m) \\
i_1{+}1 < i_2 < \cdots < i_m = n
\end{array} \\
\hline
\textbf{Relations} &
\begin{array}{l}
\rule{0pt}{3ex}
\,[s_i, s_j] = u_{ij}(s_1, \cdots, s_{i-1}) \\
\,[s_1, s_p] = 1 \\
\,1 < i < j \le n, 1 \le p \le n
\end{array}
&
\begin{array}{l}
\rule{0pt}{3ex}
\,[\overline{s_i}, \overline{s_j}] = \overline{u_{ij}(s_1, \cdots, s_{i-1})} \\
\,[\overline{s_{i_1{+}1}}, \overline{s_p}] = 1 \\
\,i_1{+}1 < i < j \le n, i_1{+}1 \le p \le n
\end{array}\\
\hline
\end{array}
\]
}

Here, each $u_{ij}(s_1, \cdots, s_{i-1})$ is a word in the generators ${s_1, \cdots, s_{i-1}}$, and we write $\overline{s_i} := s_i G_1$ and $\overline{u_{ij}} := u_{ij} G_1$. Moreover, for each $1 \le k \le m$, we have $G_k / G_{k-1} \cong \mathbb{Z}^{i_k - i_{k-1}}$, so every $g \in G$ can be uniquely written in the form $g = s_n^{a_n} s_{n-1}^{a_{n-1}} \cdots s_1^{a_1}$ for some integers $a_1, \cdots, a_n$, establishing a one-to-one correspondence between $G$ and $\mathbb{Z}^n$. Similarly, each element in $G/G_1$ has a unique normal form $\overline{s_{n}}^{x_n} \cdot \overline{s_{n-1}}^{x_{n-1}} \cdots \overline{s_{i_1{+}1}}^{x_{i_1{+}1}}$, which we fix as the standard coordinate representation in both $G$ and $G/G_1$.

Let $h\colon \mathbb{N}^+ \to \mathbb{N}^+$ be a map, and let $k_0 < k_1 < \cdots < k_m$. Since $G / G_1$ has the $(m{-}1,0)$-DTUT property, we assume that the associated constants are $(F_2, \cdots, F_m)$ and consider the subsequence $k_1 < k_2 < \cdots < k_m$. Then for each fixed $(l_2, \ldots, l_m) \in \prod_{j=2}^m \{1,\cdots, h(k_j)\}$, there exists a $k_1$-disjoint uniformly bounded family $\mathcal{C}_{G/G_1}^{l_2,\cdots,l_m}$, and moreover for each $2 \le r \le m$ and $0 \le j_r \le F_r$, a $k_r$-disjoint uniformly bounded family $\mathcal{D}_{G/G_1,r,j_r}^{l_2,\cdots,l_m}$, such that these families together cover $G/G_1$. Moreover, for any fixed $(l_2, \cdots, l_{r{-}1}, l_{r{+}1}, \cdots, l_m)$ and fixed $(r, j_r)$, the family $\bigcup_{l_r = 1}^{h(k_r)} \mathcal{D}_{G/G_1,r,j_r}^{l_2,\cdots,l_m}$ remains $k_r$-disjoint.

Following the strategy in Proposition~\ref{prop:H_3(Z) has (2,0)DTUT}, the key idea is to consider a natural coarse embedding $\phi \colon (G_{1}, \{s_1, s_2, \cdots, s_{i_{1}}\}) \to G$, and for each fixed $(l_2, \ldots, l_m)$, to construct a map $\psi^{l_2,\cdots,l_m} \colon G/G_1 \to G$ that preserves both disjointness and uniform boundedness of families in $G/G_1$.

The map $\psi^{l_2, \cdots, l_m}$ can be constructed as follows: 
for each subset $D^{l_2, \cdots, l_m}$ contained in one of the families that cover $G/G_1$ as described above, the nilpotency of $G$ and the standard coordinate representation in $G$ and $G/G_1$ ensure the existence of a correction set 
\[
\omega(G_1, D^{l_2, \cdots, l_m}) := \big\{ \omega(G_1,\, \vec{x},\, D^{l_2, \cdots, l_m}) ~\mid~ \vec{x} \in D^{l_2, \cdots, l_m} \big\} \subset G_1,
\]
where each $\omega(G_1,\, \vec{x},\, D^{l_2, \cdots, l_m})$ is a word in the generators $s_1, \cdots, s_{i_1}$, with coefficients determined by the coordinates of $\vec{x}$ and the commutator relations in $G$. For any element $\vec{x} = \overline{s_{n}}^{x_n} \cdot \overline{s_{n-1}}^{x_{n-1}} \cdots \overline{s_{i_1{+}1}}^{x_{i_1{+}1}}\in D^{l_2, \cdots, l_m}$, define
\[
\psi^{l_2, \cdots, l_m}(\vec{x}) := s_{n}^{x_n} \cdot s_{n-1}^{x_{n-1}} \cdots s_{i_1{+}1}^{x_{i_1{+}1}} \cdot \omega(G_1,\, \vec{x},\, D^{l_2, \cdots, l_m}).
\]
This establishes a ``local'' coarse embedding as required, in the sense that the prefix $s_{n}^{x_n} \cdot s_{n-1}^{x_{n-1}} \cdots s_{i_1{+}1}^{x_{i_1{+}1}}$ ensures the disjointness  (i.e., $d_G(\psi^{l_2, \cdots, l_m}(\vec{x}), \psi^{l_2, \cdots, l_m}(\vec{x}'))\ge d_{G/G_1}(\vec{x}, \vec{x}')$), while the correction term $\omega(G_1, D^{l_2, \cdots, l_m})$ preserves uniform boundedness of the above families. We assume that under $\psi^{l_2, \cdots, l_m}$, these families become $M_{G/G_1}^{\psi^{l_2, \cdots, l_m}}(k_m)$-bounded.

Since $\phi \colon (G_1, \{s_1, \cdots, s_{i_1}\}) \to G$ is a natural coarse embedding, there exists a scale $P(k_m)$ (depending on $k_m$, $\phi$, and $\psi^{l_2,\cdots,l_m}$), such that the families constructed below satisfy the required disjointness conditions. For instance, any $P(k_m)$ satisfying $\rho_-(P(k_m)) > k_m + M_{G/G_1}^{\psi^{l_2,\cdots,l_m}}(k_m)$ suffices, where $\rho_-$ is the lower control function of $\phi$. 

Since $(G_1, \{s_1, \cdots, s_{i_1}\})$ has the $(1,0)$-DTUT property with constant $F_1 = i_1 \times 2^{i_1}{-}1$, it follows that for each $1 \le l_1 \le h(k_1)$, there exists a $P(k_m)$-disjoint family $\mathcal{C}_{G_1}^{l_1}$, and for each $0 \le i \le F_1$, a $P(k_m)$-disjoint family $\mathcal{D}_{G_1,i}^{l_1}$, as required by the $(1,0)$-DTUT property.

Then by taking ``pointwise'' products between corresponding families in $G_1$ and $G/G_1$, we obtain:
\begin{itemize}
  \item one $k_0$-disjoint family: $\psi^{l_2,\cdots,l_m}\big( \mathcal{C}_{G/G_1}^{l_2,\cdots,l_m} \big) \cdot \phi\big( \mathcal{C}_{G_1}^{l_1} \big)$, 
  \item $(F_1+1)$ many $k_1$-disjoint families: $\left\{
  \psi^{l_2,\cdots,l_m}\big( \mathcal{C}_{G/G_1}^{l_2,\cdots,l_m} \big) \cdot \phi\big( \mathcal{D}_{G_1,i}^{l_1} \big)
  \right\}_{0 \le i \le F_1}$, 
  \item for each $2 \le r \le m$ and each combination of $\mathcal{C}^{l_1}_{G_1}$ or $\mathcal{D}^{l_1}_{G_1,i}$ ($0 \le i \le F_1$) with $\mathcal{D}_{G/G_1, r, j_r}^{l_2,\cdots,l_m}$ ($0 \le j_r \le F_r$), we obtain a total of $(F_1 + 2)\times (F_r+1)$ many $k_r$-disjoint families.
\end{itemize}

These families are constructed for each $1 \le l_t \le h(k_t)$ with $1 \le t \le m$. Moreover,
\begin{itemize}
  \item The family $\bigcup_{l_1=1}^{h(k_1)} \psi^{l_2,\cdots,l_m}( \mathcal{C}_{G/G_1}^{l_2,\cdots,l_m} ) \cdot \phi( \mathcal{D}_{G_1,i}^{l_1} )$ is $k_1$-disjoint for each $0 \le i \le F_1$.
  \item For each $2 \le r \le m$, both families
  \[
  \bigcup_{l_r=1}^{h(k_r)} \psi^{l_2,\cdots,l_m}( \mathcal{D}_{G/G_1, r, j_r}^{l_2,\cdots,l_m} ) \cdot \phi( \mathcal{C}_{G_1}^{l_1} ),
  \quad
  \bigcup_{l_r=1}^{h(k_r)} \psi^{l_2,\cdots,l_m}( \mathcal{D}_{G/G_1, r, j_r}^{l_2,\cdots,l_m} ) \cdot \phi( \mathcal{D}_{G_1,i}^{l_1} )
  \]
  are $k_r$-disjoint for each $0 \le j_r \le F_r$ and $0 \le i \le F_1$.
\end{itemize}

All families constructed above are uniformly bounded.
This completes the inductive step and proves that $G$ satisfies the $(m,0)$-DTUT property.
\end{proof}

Note that every finitely generated virtually nilpotent group contains a torsion-free nilpotent subgroup of finite index (\textnormal{cf.} \cite{Baumslag1971,Segal1983}). 
By Proposition~\ref{prop:(m,n)-DTUT is a coarse inv}, it follows that every finitely generated virtually nilpotent group has the DTUT property.

In 1981, Gromov proved that a finitely generated group $G$ has polynomial growth if and only if it is virtually nilpotent (\textnormal{cf.} \cite{Gromov1981}). As a consequence, we obtain the following corollary.
\begin{corollary}\label{coro:gps with poly grow has DTUT}
Every finitely generated group with polynomial growth has the $(m,0)$-DTUT property for some $m\in\mathbb{N}$. 
\end{corollary}

\subsection{Free Groups, Hyperbolic Groups}
Let $\mathbb{F}_2$ be the free group with two generators.
\begin{lemma}(\textnormal{cf.} \cite{BD2011})\label{lemma: asdim F_2=1}
asdim $\mathbb{F}_2=1$.
\end{lemma}
\begin{proof}
Take any non-negative integer $m$, and let 
\begin{align*}
{{{^m}A_k^p}}:=\big\{a\in \mathbb{F}_2 \mid l_{\mathbb{F}_2}(a)\in [(k-1)m+p,~km+p)\big\},
\end{align*}
for all $p\in\mathbb{N}$ and $k\in\mathbb{N}^+$. Then for each fixed $p \ge 0$ and $k \ge 2$, or $p \ge \frac{m}{2}$ and $k = 1$, we define an equivalence relation on ${{^m}A_k^p}$ as follows: 
define $x\sim y$ in ${{^m}A_k^p}$, if the geodesics $[1_{\mathbb{F}_2}, x]$ and $[1_{\mathbb{F}_2}, y]$ in $Cay(\mathbb{F}_2)$ contain the same point $z$ with $l_{\mathbb{F}_2}(z)\ge p+m(k-\frac{3}{2})$ (see Figure \ref{fig: the equivalence relation} for an example). It is easy to check that
\begin{itemize}
  \item "$\sim$" is an equivalence relation.
  \item the equivalence classes are $3m$-bounded.
  \item elements from distinct classes are at least $m$ apart.
\end{itemize}

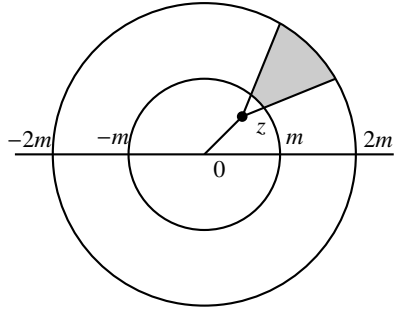
\begin{figure}[h]
\centering
\begin{tikzpicture}[scale=1, every node/.style={font=\small}]

    \tikzstyle{axis} = [thick, ->, >=stealth]
    \tikzstyle{circle_style} = [thick]
    \tikzstyle{line_style} = [thick]
    
    \draw[circle_style] (0,0) circle (1cm); 
    \draw[circle_style] (0,0) circle (2cm); 
    
    \draw[line_style] (-2.5,0) -- (2.5,0); 
    
    \draw[line_style] (0,0) -- (0.5,0.5) node[pos=1.1, below right] {$z$}; 
    \fill (0.5, 0.5) circle (2pt);

    \pgfmathsetmacro{\x}{cos(39)} 
    \pgfmathsetmacro{\y}{sin(39)} 
    \pgfmathsetmacro{\q}{cos(51)} 
    \pgfmathsetmacro{\w}{sin(51)} 
    \fill[opacity=0.2] (\x,\y) -- (1.732, 1) arc[start angle=30, end angle=60, radius=2cm]-- (1, 1.732) -- (\q,\w) arc[start angle=51, end angle=39, radius=1cm]--  cycle;
    \draw[line_style] (0.5,0.5) -- (1.732, 1);
    \draw[line_style] (0.5,0.5) -- (1, 1.732);

    \node at (1.2, 0.2) {$m$};
    \node at (-1.2, 0.2) {$-m$};
    \node at (2.3, 0.2) {$2m$};
    \node at (-2.3, 0.2) {$-2m$};
    \node at (0.2, -0.2) {$0$};

\end{tikzpicture}
\caption{Examples of the equivalence relation for $k=2$ and $p=0$.}
\label{fig: the equivalence relation}
\end{figure}

Let ${^m}\mathcal{F}_k^p$ denote the family of equivalence classes in ${^m}A_k^p$ for all $(k,p)$ satisfying either $k \ge 2$ and $p \ge 0$, or $k = 1$ and $p \ge \frac{m}{2}$.  
For the remaining case where $k = 1$ and $p < \frac{m}{2}$, define
${^m}\mathcal{F}_1^p :=\{B_{\mathbb{F}_2}(1_{\mathbb{F}_2}, p + m)\}$. Take 
\begin{align*}
{^m}\mathcal{F}_{<0>}:=\bigcup_{k\in\mathbb{N}^+,\ k\text{ odd}}{^m}\mathcal{F}_k^0,\qquad 
{^m}\mathcal{F}_{<1>}:=\bigcup_{k\in\mathbb{N}^+,\ k\text{ even}}{^m}\mathcal{F}_k^0.
\end{align*}
Then both ${^m}\mathcal{F}_{<0>}$ and ${^m}\mathcal{F}_{<1>}$ are $m$-disjoint and $3m$-bounded. It follows that $\asdim\mathbb{F}_2\le1$. Since $\mathbb{Z}<\mathbb{F}_2$, and $\asdim\mathbb{Z}=1$, we obtain that $\asdim\mathbb{F}_2=1$. 
\end{proof}

In $\mathbb{F}_2$, note that all elements lying in the annulus of word length between $pm$ and $(p+1)m$ (for $p, m \in \mathbb{N}^+$) form a $3m$-disjoint family of subsets. Therefore, similarly to the case of $\mathbb{Z}$, we conclude that $\mathbb{F}_2$ also has the $(1,0)$-DTUT property.

We use the notations established in the proof of Lemma~\ref{lemma: decomp for the infinite sum (pre for wr prod)}.

\begin{proposition}\label{prop:F_2 has (1,0)-DTUT}
The free group $\mathbb{F}_2$ has the $(1,0)$-DTUT property.    
\end{proposition}

\begin{proof}
Let $h \colon \mathbb{N}^+ \to \mathbb{N}^+$ be an arbitrary map.  
For every pair $k_0 < k_1$ with $k_0, k_1 \in \mathbb{N}^+$, every $l \in \{1,2,\cdots,h(k_1)\}$, and every $n\in \mathbb{N}^+$, define  
\begin{align*}
&\mathcal{D}^{l,n}:={^{k_1}}\mathcal{F}^0_{1+(n-1)\times 2h(k_1)+2(l-1)},\\
&\mathcal{C}^{l,n}:={^{2h(k_1)\times k_1}}\mathcal{F}^{2(l-1)k_1}_{n}\setminus \big(\bigcup\mathcal{D}^{l,n}\big).
\end{align*}
Let 
\begin{align*}
&\mathcal{D}^l:=\bigcup_{n\in \mathbb{N}^+}\mathcal{D}^{l,n},\\
&\mathcal{C}^l:=\bigg(\bigcup_{n\in \mathbb{N}^+}\mathcal{C}^{l,n}\bigg)\cup {^{2(l-1)k_1}}\mathcal{F}^0_1.
\end{align*}

We can verify that the following conditions are satisfied:
\begin{itemize}
  \item For every $1\le l\le h(k_1)$, the family $\mathcal{C}^l\cup \mathcal{D}^l$ covers $\mathbb{F}_2$. Moreover, both $\mathcal{C}^l$ and $\mathcal{D}^l$ are $k_1$-disjoint, and $(3\times 2h(k_1)\times k_1)$-bounded.
  \item The family $\bigcup_{l=1}^{h(k_1)}\mathcal{D}^l$ is $k_1$-disjoint, and $(3\times 2h(k_1)\times k_1)$-bounded.
\end{itemize}
\end{proof}

Similar to the discussion above, we observe that any tree (equipped with the standard path metric) has the $(1,0)$-DTUT property. Then, by Proposition~\ref{prop:(m,n)-DTUT is a coarse inv} and Proposition~\ref{prop:(m,n)-DTUT is closed under product}, any group that admits a coarse embedding into a finite product of trees also has the $(1,0)$-DTUT property. 

In particular, J.~Mackay and A.~Sisto showed in 2013 (\textnormal{cf.} \cite{MS2013}) that every hyperbolic group admits a quasi-isometric embedding into a finite product of metric trees. Moreover, any group that is hyperbolic relative to a collection of subgroups—each of which quasi-isometrically embeds into a finite product of metric trees—also admits a quasi-isometric embedding into such a product. 

Therefore, we obtain the following corollary.

\begin{corollary}\label{coro:hyper gps has (1,0)-DTUT}
Let $\Gamma$ be a group that coarsely embeds into a finite product of trees. Then $\Gamma$ has the $(1,0)$-DTUT property.
In particular, every hyperbolic group has the $(1,0)$-DTUT property, and so does any group that is hyperbolic relative to a collection of subgroups, each of which quasi-isometrically embeds into a finite product of trees.
\end{corollary}

\section{Questions and Perspectives}

While investigating the conditions under which a group $W \in \WR(G, H \curvearrowright I)$ admits asymptotic property C, we observe that the key requirement on $G$ lies in its possession of the disjointness of the tiling through uniform translation (DTUT) property.

In Section~\ref{sec:Groups with DTUT}, we first verified that $\mathbb{Z}^n$ satisfies the $(1,0)$-DTUT property for every $n\in\mathbb{N}^+$, and that the integral Heisenberg group satisfies the $(2,0)$-DTUT property. We then established that any finitely generated torsion-free nilpotent group of class $n$ admits the $(n,0)$-DTUT property. Since DTUT is a coarse invariant, it follows that all groups with polynomial growth possess the DTUT property. 

We further identified a special class of groups—namely, hyperbolic groups—which also satisfy the $(1,0)$-DTUT property. This is particularly interesting given that 
(non-elementary) hyperbolic groups are typically centerless, in the sense that the center of any non-elementary hyperbolic group is finite.
This leads to the following question:
\begin{question}
Under what conditions can centerless groups have the DTUT property?
\end{question}

Moreover, for amenable groups, one can construct tilings of the group using Følner sets from a Følner sequence. Specifically, in 2019, T. Downarowicz, D. Huczek, and G. Zhang proved that for any infinite countable amenable group \( G \), any \( \varepsilon > 0 \), and any finite subset \( K \subset G \), there exists a tiling of \( G \) (i.e., a partition into ``tiles'' using only finitely many ``shapes'') such that all tiles are \( (K, \varepsilon) \)-invariant (\textnormal{cf.} \cite{DHZ2019}).

In their construction of exact tilings, the shapes used in the partition are carefully chosen from a Følner sequence. These shapes are then modified through several refinement steps, including enforcing disjointness, achieving precise covering, and constructing congruent tilings with nested structure. Crucially, these modifications are quantitatively controlled so as to preserve the essential invariance properties of the original Følner sets. In particular, the deviation between each modified shape and its original Følner set is kept small in symmetric difference relative to the size of the set. 

Naturally, this leads us to consider amenable groups as potential candidates for possessing the DTUT property. This motivates the following question:

\begin{question}
Do all countable discrete amenable groups admit the DTUT property?
\end{question}


\begin{ack}
I would like to thank Y. Wu and J. Zhu for their work on establishing the asymptotic property C of $\mathbb{Z} \wr \mathbb{Z}$, whose methods inspired the generalizations explored in this paper.
\end{ack}

\bibliographystyle{emss}  
\bibliography{main}        

\end{document}